# Optimal harvesting policy for biological resources with uncertain heterogeneity for application in fisheries management

Hidekazu Yoshioka [a,*]

[a] Japan Advanced Institute of Science and Technology, 1-1 Asahidai, Nomi, Ishikawa 923-1292, Japan

* Corresponding author. E-mail: yoshih@jaist.ac.jp; Tel.: 81-761-51-1745, ORCID: 0000-0002-5293-3246

**Abstract**

Conventional harvesting problems for natural resources often assume physiological homogeneity of the body length/weight among individuals. However, such assumptions generally are not valid in real-world problems, where heterogeneity plays an essential role in the planning of biological resource harvesting. Furthermore, it is difficult to observe heterogeneity directly from the available data. This paper presents a novel optimal control framework for the cost-efficient harvesting of biological resources for application in fisheries management. The heterogeneity is incorporated into the resource dynamics, which is the population dynamics in this case, through a probability density that can be distorted from the reality. Subsequently, the distortion, which is the model uncertainty, is penalized through a divergence, leading to a non-standard dynamic differential game wherein the Hamilton–Jacobi–Bellman–Isaacs (HJBI) equation has a unique nonlinear partial differential term. Here, the existence and uniqueness results of the HJBI equation are presented along with an explicit monotone finite difference method. Finally, the proposed optimal control is applied to a harvesting problem with recreationally, economically, and ecologically important fish species using collected field data.

**Keywords**

Optimal control; heterogeneity; Hamilton–Jacobi–Bellman–Isaacs equation; finite difference method; *Plecoglossus altivelis altivelis*

**Data availability statement**

Data will be available upon reasonable request to the corresponding author.

**Funding statement**

This work was supported by the Japan Society for the Promotion of Science (grant numbers 22K14441 and 22H02456).

**Conflict of interest disclosure**

The author has no known conflict of interest.

---

Abbreviations: Hamilton–Jacobi–Bellman–Isaacs (HJBI); Hii River Fishery Cooperative (HRFC)

**Ethics approval statemen**

N.A.

**Patient consent statement**

N.A.

**Permission to reproduce material from other sources**

N.A.

**Clinical trial registration**

N.A.

**Declaration of generative AI in scientific writing**

The authors did not use any generative AI technologies for writing this manuscript.

**Author contributions**

Only the corresponding author contributed to this work.

**Acknowledgements**

The authors thank the Hii River Fishery Cooperative for providing the data of the body weights of the fish *Plecoglossus altivelis altivelis* in the Hii River, Japan. Valuable comments from the reviewers improved the manuscript.

# 1. Introduction

## 1.1 Research background

Biological resources, such as fishes, crops, and livestock, are essential for the sustainability of human lives [1–4]. The aim of biological resource management is to maximize the utility while minimizing the harvesting cost [5,6] as well as the disutility that is triggered by resource exhaustion or extinction [7,8]. Thus, optimal control theory provides normative descriptions to deal with the general resource management problem (e.g., [9,10]).

Many classical optimal control problems for biological resource harvesting do not account for physiological heterogeneity, such as the timber diameter, body length, and body weight, among individuals, possibly to ensure modeling simplicity. However, such heterogenous qualities or individual differences are not negligible in real-world problems. In the simplest case, individual differences emerge as a probability distribution with respect to a physiological variable such as the body weight [11] and the distance between circuli in a scale [12]. Ghwila and Willms [13] modeled physiologically structured population dynamics and the stability thereof under prescribed harvesting policies.

The emergence of heterogeneity in biological resource dynamics mathematically involves another dimension that is often infinite-dimensional, which increases the model complexity, although this aspect can be resolved analytically if the problem at hand admits a tractable structure. Examples include, but are not limited to, the graphon linear-quadratic game whose resolution is reduced to the Riccati equation [14], and the moment approximation closure [15]. However, most problems that are encountered in applications do not admit such a useful mathematical structure. Furthermore, the engineering applications of the aforementioned models have not been sufficiently explored in the literature. Therefore, a mathematical model for biological resource management is required that accounts for heterogeneity while allowing the associated control problem to be solved efficiently, without resorting to complex numerical methods.

Certain control problems are essentially heterogenous in space [16,17] or age [18,19]. However, our focus is on the physiological heterogeneity among individuals in a habitat, as considered in some lumped fisheries management problems in a water body such as an aquaculture system or a river reach. Moreover, social heterogeneity [20] is an important concept for explaining and controlling specific population dynamics such as disease spreading. Ensemble control can also potentially handle heterogenous populations when the total number of individuals is not large [21].

The modeling of biological resource dynamics, regardless of whether they are homogenous or heterogenous, in applications faces the limitation of data availability [22,23]. This data limitation creates a bottleneck in decision-making in resource management [24, 25]. The effects of model uncertainty can be evaluated using a Monte-Carlo simulation, as has been done in engineering research [26,27]; however, this approach is not easily implementable in control problems that are subject to model uncertainty even under a homogenous case, because it is inefficient for solving the control problem numerically a vast number of times.

Multiplier robust control [28] is a differential game formalism of optimal control problems that is subject to model uncertainty, and has mainly been applied to stochastic control problems in which the underlying probability measure is distorted owing to insufficient knowledge of the decision-makers [29–31]. Barnett et al. [32] parameterized the model uncertainty using prior probability weights, which can be considered as a form of uncertain heterogeneity in harmony with the control problem. A theoretical as well as practical advantage of the multiplier robust approach is that the model uncertainty can be represented by suitably perturbing the objective function and system dynamics of the targeted control problem within the context of the major optimality principle, such as the dynamic programming and maximum principle. Nevertheless, to the best of the authors' knowledge, this approach has not been incorporated into biological resource management with physiological heterogeneity.

### 1.2 Research objective and contribution

The objectives of this study are twofold: the formulation of a finite-horizon optimal control problem for biological resource physiological heterogeneity, and the application thereof to a real case. To achieve the first objective, we consider a harvesting problem of a population in a habitat in which the individuals are physiologically heterogenous, with an emphasis on fisheries management. The proposed model couples the population dynamics coupled with a continuum of growth curves representing the physiological heterogeneity, leading to a unique hybrid system that has not been considered. More specifically, the heterogeneity is represented using a probabilistic parameterization of the logistic-type growth model, where the growth rate and/or maximum body weight are randomly distributed among the individuals (e.g., [33,34]). The uncertainty here concerns the heterogeneity: the probability distribution of the physiological heterogeneity. Therefore, in the proposed model there is no uncertainty if there is no heterogeneity, which simplifies the theoretical formulation but in a nontrivial way as discussed below.

Note that the approaches in the literature [20, 21] are based on microscopic individual-based models, whereas ours is an aggregated one that deals with the population as a whole. The two approaches are therefore fundamentally different with each other. The approach based on the incomplete information [48] assumes that a model parameter is possibly better estimated as the time is elapsed due to a statistical filter, while ours assume that this mechanism is absent. Instead, we assume in this study that the decision-maker can manage the fishery resource considering model uncertainties through the uncertainty aversion term in the objective function.

The novelty of this study lies in the combination of multiplier robust control with this randomized formalism in a control problem to account for possible misspecifications of the probability measure when generating the randomness. The degree of randomness is represented by the Kullback–Leibler divergence, which is a widely used statistical divergence (e.g., [35]). We demonstrate that solving this control problem formally can be reduced to determining a proper solution to the Hamilton–Jacobi–Bellman–Isaacs (HJBI) equation. In this equation, the consideration of the model uncertainty emerges as a unique nonlinear term that potentially raises issues in both theory and computation. We demonstrate that, despite its complicated form, the HJBI equation admits a Lipschitz continuous Hamiltonian, to which the comparison argument to prove the unique solvability in a viscosity sense (e.g., [36,37]) and a convergent finite difference discretization (e.g., [38,39]) can be applied with modifications.

Our HJBI equation is related to problems that arise in machine learning, especially in reinforcement learning, where the decision-maker should select randomized controls (e.g., [40]). Exploratory control [41], maximum entropy optimal control [42], and feedback relaxed control [43] lead to specific nonlinear and nonlocal optimality equations corresponding to the randomized control. Our HJBI equation differs because the control problem is aimed at maximizing some utility under uncertain randomness in the population dynamics; hence, the system dynamics, rather than the control variable, is randomized. In this sense, the linear-quadratic regulator of the uncertain system matrix [44] is closer to ours, although this methodology depends on the maximum principle rather than dynamic programming. Nevertheless, our and the aforementioned optimality equations share some similarities, such as the viscosity property, which facilitate the understanding of all their formalisms. It should be noted that, to the best of our knowledge, our control problem has not yet been studied, even with the absence of model uncertainty. Furthermore, our control problem is not expected to satisfy the so-called Isaacs condition (i.e., the order of maximization and minimization can be changed in a differential game) owing to the lack of convexity of the Hamiltonian.

We also present a simple finite difference method for discretizing the HJBI equation to provide demonstrative examples using real data. Finally, we apply the proposed control problem to a harvesting problem of the inland fishery resource Ayu *Plecoglossus altivelis altivelis*, which is one of the most recreationally, economically, and ecologically important fish species in Japan [45,46]. Although the authors have studied harvesting problems of this fish, their models have assumed only homogeneous populations [47,48]. The proposed mathematical framework that accounts for physiological heterogeneity extends these methodologies. The population dynamics of the fish is estimated from the unique collected data of the body weights of the individuals since 2020. The computational results demonstrate the optimal harvesting policy of *P. altivelis* depending on the control objective as well as the potential model uncertainty. We demonstrate that the uncertainty acting on the heterogeneity gives the worst-case underestimation of the probability density of the body weights of the individual fishes in the population according to the uncertainty aversion of the decision-maker. Different probability densities having different skewness are examined to see their impacts on the harvesting policy and the worst-case growth curve. This work covers the formulation, analysis, and application of a novel optimal control problem.

The remainder of this paper is organized as follows: The mathematical framework that is employed in this study is introduced in **Section 2**. The mathematical analysis with a focus on the HJBI equation is presented in **Section 3**, along with the finite difference method and its theoretical analysis results. **Section 4** outlines the model application to real data, focusing on the harvesting problem of *P. altivelis*. A summary and future perspectives of our research are presented in **Section 5**. Proofs of the propositions and lemmas are provided in **Appendix A**. Few additional computational results are presented in **Appendix B**.

## 2. Mathematical model

We consider a harvesting problem of a biological resource in a closed habitat in a finite horizon during which spawning (i.e., a population increase) does not occur, which is effectively considered as the terminal utility in the objective function.

### 2.1 Resource dynamics
#### 2.1.1 Population dynamics

This study focuses on fisheries management, and the utility is evaluated through the biomass of aggregated individuals. The resource dynamics in this paper contains two continuous-time variables: the population (i.e., total number of individuals) and a probabilistically distributed body weight. The population $(N_t)_{t \geq 0}$ evolves according to the differential equation

$$\frac{dN_t}{dt} = -R(N_t)N_t - c_t, \quad 0 < t < \tau_{0,n}, \quad N_0 = n > 0, \tag{1}$$

where $t$ is the time, $R : [0, +\infty) \to [0, +\infty)$ is the mortality rate that is bounded and positive, increasing, and Lipschitz continuous: $|R(n_1) - R(n_2)| \leq L_R |n_1 - n_2|$ for all $n_1, n_2 \geq 0$ with a constant $L_R > 0$, $(c_t)_{t \geq 0}$ represents the harvesting rate that belongs to the admissible set

$$\mathbb{A} = \left\{ (c_t)_{t \geq 0} \,\middle|\, c_t \text{ is measurable and } 0 \leq c_t \leq \bar{c} \ (t \geq 0) \right\} \tag{2}$$

with a constant $\bar{c} > 0$, and $\tau_{t,n}$ is the stopping time such that

$$\tau_{t,n} = \inf \left\{ \tau \geq t \,\middle|\, N_\tau = 0 \text{ and } N_t = n \right\}. \tag{3}$$

We set $N_s = 0$ for all $s \geq \tau_{t,n}$. The population is approximated to be a real variable rather than an integer, which can be reasonable if it is significantly larger than $O(10^0)$, such as $O(10^4)$, as considered in our application. The population-dependent mortality rate covers both constant and state-dependent cases, with the latter able to consider the density-dependent mortality if necessary (e.g., [49,50]). The constant $\bar{c}$ that serves as the upper bound of the control variable $c$ represents a technological constraint. It also ensures the unique existence of the solution to (1), as indicated below. The set $\mathbb{A}$ is non-empty as the constant control $c \equiv 0$ belongs to it.

The differential equation (1) is solved for $0 < t < \tau_{0,n}$ as follows:

$$N_t = n \exp\left(-\int_0^t R(N_s) ds\right) - \int_0^t \exp\left(-\int_s^t R(N_u) du\right) c_s ds, \quad 0 < t < \tau_{0,n}, \tag{4}$$

which implicitly defines the unique continuous solution to (1) because the right-hand side of (4) is well defined and $R$ is bounded and Lipschitz continuous (e.g., Theorem 5.5 in Chapter 3, Section 5 of [51]; see also conditions (A$_0$), (A$_1$), and (A$_3$) of this work). The population $N$ is set to 0 after the hitting time $\tau_{0,n}$. In this manner, $N_t$ is a continuous decreasing function in $[0, +\infty)$, where the decreasing property is immediate from the non-positivity of the right-hand side of (1).

### 2.1.2 Growth dynamics

Another component of the resource dynamics is the continuum of growth equations, which is simply referred to as the growth dynamics. We assume the widely-used logistic model (e.g., [34]):

$$X_t(u) = K(u)\left[1 + \left(\frac{K(u)}{x} - 1\right)\exp(-r(u)t)\right]^{-1}, \quad t \geq 0, \quad (5)$$

where $u \in [0,1]$ parameterizes a continuum of individuals with a distributed growth rate $r(u)$ and maximum body weight $K(u)$ such that $r, K$ are bounded and positive, and strictly increasing in $[0,1]$. Furthermore, the initial condition $x > 0$ is sufficiently small such that $x \leq K(0)$. For each $u \in [0,1]$, the right-hand side of (5) exists as a uniformly continuous function of time $t \geq 0$. We assume that there exists a probability density function $p(u)$ such that the averaged weight in the population is defined as

$$\bar{X}_t = \int_0^1 X_t(u) p(u) \mathrm{d}u, \quad t \geq 0. \quad (6)$$

Subsequently, the total biomass of the population at time $t \geq 0$ is given by $N_t \bar{X}_t$. Similarly, the unit-time harvested biomass is given by $c_t \bar{X}_t$.

***Remark 1*** Employing the logistic growth curve is for analytical simplicity. Other growth curves, such as the Gompertz and generalized logistic curves (e.g., [34]), can be used instead of (5). Theoretically, $X_t(u)$ is required to be uniformly continuous and strictly increasing with respect to $t \geq 0$ as well as strictly positive and strictly bounded from above for all $t \geq 0$ and $u \in [0,1]$; that is, we require $0 < \tilde{X}_1 \leq \min_{u \in [0,1]} X_0(u) < \max_{u \in [0,1]} \lim_{t \to +\infty} X_t(u) \leq \tilde{X}_2 < +\infty$ with constants $\tilde{X}_1, \tilde{X}_2 > 0$. A key point in using the S-shaped curve is its profile having an inflection point, with which the model can be reasonably fitted against data as demonstrated in this study. Another key point, which will be critical in applications, is its boundedness such that the growth curve eventually saturates as the time is elapsed. Growth curves that are simpler than the S-shaped one, such as an exponentially increasing curve, will be easier to identify and operate, whereas its unboundedness will overestimate the growth curve and the biomass. Similarly, one will use a variety of distributions $p$. In **Section 4**, we argue that the simplest distribution will be the uniform, one and beta distributions can be used for unimodal alternatives; multi-model ones can also be used if necessary, but is not the case for the application to *P. altivelis*.

### 2.2 Model uncertainty

The model uncertainty is represented by a Radon–Nikodym derivative as a positive and measurable continuous-time field $(\phi_t(\cdot))_{t \geq 0}$, which is referred to as the uncertainty field as well in this paper, such

that the integrability condition is satisfied:

$$\int_0^1 \phi_t(u) p(u) du = 1, \quad t \geq 0 \tag{7}$$

and the Kullback–Leibler divergence $\mathbb{D}$ (e.g., [35]) exists:

$$\mathbb{D}(\phi_t) = \int_0^1 (\phi_t(u) \ln \phi_t(u) - \phi_t(u) + 1) p(u) du \in [0, +\infty), \quad t \geq 0. \tag{8}$$

The uncertainty field can be considered as an uncertainty that is induced by incomplete knowledge of the growth dynamics. For later use, we set the admissible set of $\phi$ as follows:

$$\mathbb{B} = \left\{ (\phi_t(\cdot))_{t \geq 0} \mid \phi_t(\cdot) \ (t \geq 0) \text{ is measurable, non-negative, and satisfies (7)-(8)} \right\}. \tag{9}$$

Note that the function $y \ln y - y + 1$ ($y \geq 0$) is strictly convex, non-negative, and globally minimized at $y = 1$ with the minimum value 0. Here, we understand $0 \ln 0 = 0$. The set $\mathbb{B}$ is non-empty as the constant field $\phi \equiv 1$ belongs to it.

According to the formulation above, the heterogeneity is modulated by the uncertainty. This implies that the uncertainty not only changes the profile of $p$ but also its moments, which further affect the biomass. This means that we focus on the uncertainty in the growth curve that eventually affects the biomass of the target fish and hence its resource management.

## 2.3 Control objective

We consider a finite-horizon control problem during a fixed period $[t, T]$ for each $0 \leq t \leq T$ with a prescribed terminal time $T > 0$. For later use, for any $t \geq 0$, $c \in \mathbb{A}$, and $\phi \in \mathbb{B}$, the population that is conditioned on $N_t = n$ at time $t$ and evolves according to differential equation (1) is denoted as $\left(N_s^{(t,n,c,\phi)}\right)_{s \geq t}$.

The objective of the control problem consists of a cumulative utility by harvesting, a cumulative harvesting cost, a terminal utility, and a cumulative penalization of the model uncertainty; for each $0 \leq t \leq T$, $n \geq 0$, $c \in \mathbb{A}$, and $\phi \in \mathbb{B}$, we set

$$J(t,n,c,\phi) := \int_t^T \left( \underbrace{\alpha (c_s \overline{X}_{s,\phi})^\beta}_{\text{Harvesting utility}} \underbrace{- \gamma c_s}_{\text{Harvesting cost}} \right) ds + \underbrace{h\left(N_T^{(t,n,c,\phi)}\right)}_{\text{Terminal utility}} + \underbrace{\mu \int_t^T \mathbb{D}(\phi_s) ds}_{\text{Penalization of uncertainty}}, \tag{10}$$

where

$$\overline{X}_t^{(\phi)} = \int_0^1 X_t(u) \phi_t(u) p(u) du. \tag{11}$$

In the first term of (10), $\alpha > 0$ is the weighting factor of the harvesting utility and $\beta \in (0,1)$ is the power index of the utility. The concavity of this term with respect to $c_s$ combined with the linear harvesting cost results in a balanced optimal harvesting rate, as demonstrated later. Moreover, it contains a tractable case $\beta = 1/2$, in which the Hamiltonian of the HJBI equation is determined analytically. In the second term, $\gamma > 0$ is the weighting factor of the harvesting cost, where the cost is proportional to the harvesting rate, which assumes that the labor that is required to harvest the fish increases linearly with respect to the population to be harvested because harvesting a larger population requires a larger effort. This linearity is another key to obtain the optimal control in a closed form. In the third term, $h:[0,+\infty) \to [0,+\infty)$ with $h(0) = 0$ is a bounded, increasing, and Lipschitz continuous function that represents the terminal utility, which would effectively suppress the resource extinction [47]. In **Section 4**, we computationally demonstrate that the use of a simple form of $h$ realizes a positive terminal population that potentially contributes to the spawning, and hence, maintains the life cycle of the biological resource. Finally, in the last term, $\mu > 0$ is a weighting factor that serves as an uncertainty aversion parameter such that a larger (resp., smaller) $\mu$ more weakly (resp., more strongly) allows for the existence of the model uncertainty, and hence, supposes a model closer to (resp., farther from) the benchmark model with $\phi \equiv 1$. Expectation (11) is a distorted counterpart of (6) according to the uncertainty field $\phi$.

We set the domain $\Omega = (0,T) \times (0,+\infty)$ and its closure $\bar\Omega = [0,T] \times [0,+\infty) \to \mathbb{R}$. Based on the objective (10), we define the value function $\Phi: \bar\Omega \to \mathbb{R}$:

$$\Phi(t,n) = \inf_{\phi \in \mathbb{B}} \sup_{c \in \mathbb{A}} J(t,n,c,\phi). \tag{12}$$

This value function is understood as a pessimistic maximum value of the control objective (10) because the right-hand side is a minimum of a maximum. According to (10), $\Phi$ should satisfy the following boundary and terminal conditions:

$$\Phi(t,0) = 0, \quad 0 < t < T \quad \text{(no utility/cost if there is no population)} \tag{13}$$

and

$$\Phi(T,n) = h(n), \quad n \geq 0 \quad \text{(the terminal value is attained at time } T\text{)}, \tag{14}$$

respectively. The boundary and terminal conditions are compatible at $(t,n) = (T,0)$ due to $h(0) = 0$.

In the following, we focus on the aforementioned case $\beta = 1/2$ as it is the most tractable case among $\beta \in (0,1)$, where the HJBI equation becomes the most tractable while its characteristics will not be critically lost.

The value function $\Phi$ is non-negative and bounded. Indeed, its non-negativity, which is a lower bound, follows by considering $J$ with constant controls in combination with the non-negativity of $\mathbb{D}$, $h$:

$$\Phi(t,n) = \inf_{\phi \in \mathbb{B}} \sup_{c \in \mathbb{A}} J(t,n,c,\phi) \geq \inf_{\phi \in \mathbb{B}} J(t,n,c,\phi)\Big|_{c \equiv 0} = \inf_{\phi \in \mathbb{B}} \mu \int_t^T \mathbb{D}(\phi_s) \mathrm{d}s = 0. \tag{15}$$

The other bound follows from

$$\begin{aligned}
\Phi(t,n) &= \inf_{\phi \in \mathbb{B}} \sup_{c \in \mathbb{A}} J(t,n,c,\phi) \\
&\leq \sup_{c \in \mathbb{A}} \left\{ \int_t^T \left( \alpha \left( c_s \overline{X}_{s,\phi} \right)^\beta - \gamma c_s \right) ds + h\left( N_T^{(t,n,c,\phi)} \right) + \mu \int_t^T \mathbb{D}(\phi_s) ds \right\} \bigg|_{\phi \equiv 1} \\
&\leq \sup_{c \in \mathbb{A}} \left\{ \int_t^T \alpha \left( c_s \overline{X}_{s,\phi} \right)^\beta ds + h(n) \right\} \bigg|_{\phi \equiv 1} \\
&\leq \sup_{c \in \mathbb{A}} \left\{ \int_t^T \alpha \left( \overline{c} \overline{X}_{s,\phi} \right)^\beta ds + h(n) \right\} \bigg|_{\phi \equiv 1} \\
&\leq T \alpha \left( \overline{c} K(1) \right)^\beta + \max_{n \geq 0} h(n) \\
&< +\infty
\end{aligned} \qquad (16)$$

## 2.4 HJBI equation

### 2.4.1 Formulation

The dynamic programming principle heuristically leads to the HJBI equation. That is, the optimality equation that governs the value function $\Phi = \Phi(t,n)$ is expressed in domain $\Omega$ as follows:

$$\frac{\partial \Phi}{\partial t} - R(n) n \frac{\partial \Phi}{\partial n} + \inf_{\phi(\cdot) > 0, \int_0^1 \phi(u) p(u) = 1} \sup_{0 \leq c \leq \overline{c}} \left\{ \begin{array}{l} -c \dfrac{\partial \Phi}{\partial n} + \alpha \left( c \int_0^1 X_t(u) \phi(u) p(u) du \right)^{1/2} \\ -\gamma c + \mu \int_0^1 \left( \phi(u) \ln \phi(u) - \phi(u) + 1 \right) p(u) du \end{array} \right\} = 0. \qquad (17)$$

This HJBI equation is subject to boundary condition (13) and terminal condition (14). We can more compactly rewrite (17) by explicitly calculating the "inf sup" part, while it should also be noted that the so-called Isaacs condition (i.e., the order of "inf" and "sup" can be exchanged in (17)) is not expected, as discussed in the following sub-section.

### 2.4.2 Isaacs condition and heuristic optimal controls

The discussion in this subsection applies to the case $\beta \neq 1/2$ as well. The Isaacs condition for (17) is not expected; that is, the equality for any $t \in (0,T)$ and $z \in \mathbb{R}$,

$$\sup_{0 \leq c \leq \overline{c}} \inf_{\phi(\cdot) > 0, \int_0^1 \phi(u) p(u) = 1} \left\{ \begin{array}{l} -c(\gamma + z) + \alpha \left( c \int_0^1 X_t(u) \phi(u) p(u) du \right)^{1/2} \\ + \mu \int_0^1 \left( \phi(u) \ln \phi(u) - \phi(u) + 1 \right) p(u) du \end{array} \right\}$$
$$= \inf_{\phi(\cdot) > 0, \int_0^1 \phi(u) p(u) = 1} \sup_{0 \leq c \leq \overline{c}} \left\{ \begin{array}{l} -c(\gamma + z) + \alpha \left( c \int_0^1 X_t(u) \phi(u) p(u) du \right)^{1/2} \\ + \mu \int_0^1 \left( \phi(u) \ln \phi(u) - \phi(u) + 1 \right) p(u) du \end{array} \right\}, \qquad (18)$$

is not expected in general. The primary reason for the failure of the Isaacs condition is the left-hand side of (18). With respect to $\phi$, the first term of the functional

$$\alpha \left( c \int_0^1 X_t(u) \phi(u) p(u) du \right)^{1/2} + \mu \int_0^1 \left( \phi(u) \ln \phi(u) - \phi(u) + 1 \right) p(u) du \qquad (19)$$

is concave, whereas the second term is convex, which suggests that the terms inside "{ }" of the left-hand side of (18) are not necessarily convex particularly for small $\mu > 0$. This implies that the inner maximization on the left-hand side of (18) is not always achieved. In contrast, for each fixed $\phi(\cdot) > 0, \int_0^1 \phi(u) p(u) = 1$, the function

$$-c(\gamma+z) + \alpha \left( c \int_0^1 X_t(u) \phi(u) p(u) \mathrm{d}u \right)^{1/2} \quad \text{for} \quad c \in [0, \bar{c}] \tag{20}$$

is concave and is maximized by

$$\hat{c}(t, z, \phi) := \max\left\{ 0, \min\left\{ \frac{\alpha^2}{4(\gamma+z)^2} \left( \int_0^1 X_t(u) \phi(u) p(u) \mathrm{d}u \right), \bar{c} \right\} \right\}. \tag{21}$$

In the next step, it is necessary to minimize

$$\alpha \left( \hat{c}(t, z, \phi) \int_0^1 X_t(u) \phi(u) p(u) \mathrm{d}u \right)^{1/2} + \mu \int_0^1 \left( \phi(u) \ln \phi(u) - \phi(u) + 1 \right) p(u) \mathrm{d}u, \tag{22}$$

which turns out to be a strictly convex optimization problem under a suitable assumption. For example, this condition is satisfied with a sufficiently large $\bar{c}$ (the technological limitation does not restrict the optimal harvesting policy) and an increasing property of $\Phi = \Phi(t, n)$ with respect to $n \geq 0$ (the value function, namely the net utility, is larger for a larger population). We present the following **Assumption 1** ($\bar{c}$ is sufficiently large) and **Lemma 1** in preparation for the more detailed discussion in **Section 3**. **Assumption 1** is assumed throughout this paper to simplify the problem at hand.

*Assumption 1*

$$\frac{\alpha^2}{4\gamma^2} K(1) \leq \bar{c}. \tag{23}$$

*Lemma 1*

*For $0 < t < T$, $z \geq 0$, and $\phi(\cdot) > 0, \int_0^1 \phi(u) p(u) = 1$, it follows that*

$$\hat{c}(t, z, \phi) = \frac{\alpha^2}{4(\gamma+z)^2} \int_0^1 X_t(u) \phi(u) p(u) \mathrm{d}u. \tag{24}$$

From **Lemma 1**, for each $z \geq 0$, we obtain

$$\inf_{\phi(\cdot)>0, \int_0^1 \phi(u) p(u)=1} \sup_{0 \leq c \leq \bar{c}} \left\{ \begin{array}{l} -c(\gamma+z) + \alpha \left( c \int_0^1 X_t(u) \phi(u) p(u) \mathrm{d}u \right)^{1/2} \\ + \mu \int_0^1 \left( \phi(u) \ln \phi(u) - \phi(u) + 1 \right) p(u) \mathrm{d}u \end{array} \right\}$$

$$= \inf_{\phi(\cdot)>0, \int_0^1 \phi(u) p(u)=1} \left\{ \frac{\alpha^2}{4(\gamma+z)} \int_0^1 X_t(u) \phi(u) p(u) \mathrm{d}u + \mu \int_0^1 \left( \phi(u) \ln \phi(u) - \phi(u) + 1 \right) p(u) \mathrm{d}u \right\}. \tag{25}$$

The second line of (25) is a convex optimization problem (e.g., Section 2 of [52]), and the infimum is

achieved by the minimizer

$$\hat{\phi}(t,z,u) = \left(\int_0^1 e^{-\frac{\alpha^2}{4\mu(\gamma+z)}X_t(u)} p(u)\mathrm{d}u\right)^{-1} e^{-\frac{\alpha^2}{4\mu(\gamma+z)}X_t(u)}. \tag{26}$$

Hence, the right-hand side of (25) becomes

$$H(t,z) := -\mu \ln\left(\int_0^1 e^{-\frac{\alpha^2}{4\mu(\gamma+z)}X_t(u)} p(u)\mathrm{d}u\right). \tag{27}$$

Subsequently, we obtain the candidate for the optimal control $(c^*, \phi^*) \in \mathbb{A} \times \mathbb{B}$ to yield the value function $\Phi$:

$$\phi_t^*(u) = \hat{\phi}\left(t, \frac{\partial \Phi}{\partial n}(t, N_t), u\right), \quad 0 < t < T, \quad 0 \leq u \leq 1 \tag{28}$$

and

$$c_t^* = \frac{\alpha^2}{4\left(\gamma + \frac{\partial \Phi}{\partial n}(t, N_t)\right)^2} \int_0^1 X_t(u) \hat{\phi}\left(t, \frac{\partial \Phi}{\partial n}(t, N_t), u\right) p(u)\mathrm{d}u, \quad 0 < t < T. \tag{29}$$

In this context, the determination of the optimal controls is reduced to the resolution of the HJBI equation. Based on (28), the pessimistic estimate of the averaged body weight is obtained as follows:

$$\bar{X}_t^{(\phi^*)} = \int_0^1 X_t(u) \phi_t^*(u) p(u)\mathrm{d}u. \tag{30}$$

Formulations (28)–(30) imply that the partial derivative $\frac{\partial \Phi}{\partial n}$ plays a crucial role in determining the optimal harvesting policy as well as the model uncertainty. Phenomenologically, a larger $\frac{\partial \Phi}{\partial n}$ leads to a smaller uncertainty aversion, and hence, a more optimistic result.

***Remark 2*** The loss of the Isaacs condition discussed above is considered due to the concavity of its first term in (19). Using a convex alternative will be theoretically possible, but a problem in our context is that the optimal control may not be obtained explicitly as in (21), which makes the model be complicated. It implies that if one uses a convex alternative utility, then the Isaacs condition can be satisfied while the model complexity would significantly increase. Moreover, in such a case, the computational algorithm used in this study will have to be significantly customized, which is beyond its scope. By contrast, we consider that using an alternative distribution $p$ (other than the uniform and beta ones examined later) does not affect the Isaacs condition because this condition is with respect to $\phi$.

***Remark 3*** Eq. (26) demonstrates that, for any probability density $p$ of the heterogeneity, the worst-case distortion of the heterogeneity is proportional to the exponential $e^{-\frac{\alpha^2}{4\mu(\gamma+z)}X_t(u)}$ given in Eq. (5). As

$X_t(u)$ is increasing in $u \in [0,1]$, this exponential is monotone and more specifically decreasing in $u$. The ratio of individuals possibly having smaller (resp., larger) body weight is therefore underestimated in the parametric way. **Figure 7** in **Section 4** and **Figure B2** in **Appendix B** visualize quantitative magnitudes of the worst-case distortion on the growth curve.

## 3. Mathematical analysis

### 3.1 Value function

We analyze the monotonicity and continuity properties of the value function $\Phi$. We first present a proposition stating that $\Phi$ is increasing in the second argument. This implies that we should focus on solutions to the HJBI equation complying with $\frac{\partial \Phi}{\partial n} \geq 0$ in some sense.

*Proposition 1*

*For each $t \in [0,T]$ and $n_1, n_2 \in [0,+\infty)$ with $n_2 \geq n_1$, it follows that*

$$\Phi(t, n_1) \leq \Phi(t, n_2). \tag{31}$$

*Remark 4* One may proceed to a verification argument of viscosity solutions, which potentially yields an existence proof of the viscosity solutions to the associated HJBI equation. We do not take this direction in this study owing of the possible lack of the Isaacs condition in our HJBI equation, which makes the problem more difficult than one that complies with this condition. Instead, we construct a continuous viscosity solution to our "modified" HJBI equation with a truncated Hamiltonian, as presented later, based on the finite difference method. Owing to **Proposition 1**, this truncation can be removed and the unique viscosity solution to the modified equation is also a viscosity solution to the HJBI equation. Later, we numerically imply that the value function is continuous but not Lipschitz continuous. Possibility of formulating an analogous problem satisfying the Isaac's condition is an open issue.

### 3.2 HJBI equation

According to **Section 2.4.2**, formally, if $\frac{\partial \Phi}{\partial n} \geq 0$, the HJBI equation (17) is rewritten as

$$\frac{\partial \Phi}{\partial t} - R(n) n \frac{\partial \Phi}{\partial n} + H\left(t, \frac{\partial \Phi}{\partial n}\right) = 0. \tag{32}$$

We wish to discuss the uniqueness of the (viscosity) solutions to the HJBI equation, which is critically dependent on the regularity of the Hamiltonian $H$. We present the following **Lemma 2**, which states that $H$ is Lipschitz continuous with respect to $z \geq 0$ and uniformly continuous with respect to $t \in [0,T]$.

*Lemma 2*

The Hamiltonian $H$ in (27) is uniformly continuous with respect to $t \in [0,T]$ for each $z \geq 0$. Furthermore, it satisfies the Lipschitz continuity

$$|H(t,z_1) - H(t,z_2)| \leq \frac{\alpha^2}{4\gamma^2} K(1) |z_1 - z_2| \quad (z_1, z_2 \geq 0), \tag{33}$$

and $\frac{\partial H}{\partial z}(t,z)$ is uniformly continuous with respect to $0 \leq t \leq T$ for each $z \geq 0$. Moreover, $H$ is decreasing in the second argument.

*Remark 5* **Lemma 2** states that our Hamiltonian admits certain regularity properties despite its apparent complexity. Another key finding is that the existence of the harvesting cost, namely the positivity of $\gamma$, is essential for its regularity. Indeed, the proof of **Lemma 2** suggests that the Lipschitz constant of (33) is sharp and becomes arbitrary large as $\gamma$ approaches $0$. This is in contrast to [53], in which a costless formulation was considered, and a simpler and more tractable objective function was used.

The uniqueness is proven through the comparison argument owing to **Lemma 2**. In combination with the increasing property of the value function $\Phi$ (**Proposition 1**), if there exists a unique viscosity solution to the modified HJBI equation

$$\frac{\partial \Phi}{\partial t} - R(n)n\frac{\partial \Phi}{\partial n} + \tilde{H}\left(t, \frac{\partial \Phi}{\partial n}\right) = 0 \tag{34}$$

that is subject to the same boundary and terminal conditions as the original one, with

$$\tilde{H}(t,z) = H(t, \max\{z,0\}), \quad 0 \leq t \leq T, \quad z \in \mathbb{R}, \tag{35}$$

this solution also solves the original HJBI equation (32) if $\frac{\partial \Phi}{\partial n} \geq 0$. More specifically, the modified Hamiltonian (35) satisfies the same Lipschitz continuity (33) for any $z_1, z_2 \in \mathbb{R}$ because the function $\max\{z,0\}$ ($z \in \mathbb{R}$) has the Lipschitz constant 1. The modified Hamiltonian is defined over $z \in \mathbb{R}$ and covers a wider range than the original one (27). Suppose that the modified equation (34) admits a unique viscosity solution. Under this assumption, if a viscosity solution to the HJBI equation (32) that satisfies $\frac{\partial \Phi}{\partial n} \geq 0$ in some sense can be found, this solution will be the unique viscosity solution to the modified equation (34) as well. This implies that the solution to the modified equation (34) satisfies $\frac{\partial \Phi}{\partial n} \geq 0$, which suggests that the modified equation (34) serves as the optimality equation of our control problem as well. This type of use of a modified equation was employed in [54] to deal with a viscosity solution with a certain monotonicity property.

We define the viscosity solutions to the HJBI equation (32). The collection of all continuous functions that are defined in a domain $D$ is denoted as $C(D)$. Similarly, the collection of all upper-semicontinuous (resp., lower-semicontinuous) functions that are defined in a domain $D$ is denoted as $USC(D)$ (resp., $LSC(D)$). The linear growth speed of the viscosity solutions can be justified owing to the boundedness (16).

## Definition 1
*(a) Viscosity super-solution*

*Let $\bar{\Phi} \in LSC(\bar{\Omega})$, which grows at most linearly with respect to the second argument. We refer to $\bar{\Phi}$ as a viscosity super-solution of (34) if, for any $(\hat{t}, \hat{n}) \in \bar{\Omega}$ and $\varphi \in C^1(\bar{\Omega})$ such that $\bar{\Phi}(\hat{t}, \hat{n}) - \varphi(\hat{t}, \hat{n}) = \min_{(t,n) \in \bar{\Omega}} \{\bar{\Phi}(t,n) - \varphi(t,n)\}$, it follows that*

$$\frac{\partial \varphi}{\partial t}(\hat{t}, \hat{n}) - R(\hat{n})\hat{n}\frac{\partial \varphi}{\partial n}(\hat{t}, \hat{n}) + \tilde{H}\left(\hat{t}, \frac{\partial \varphi}{\partial n}(\hat{t}, \hat{n})\right) \leq 0 \quad \text{when} \quad (\hat{t}, \hat{n}) \in \Omega, \tag{36}$$

$\bar{\Phi}(\hat{t}, \hat{n}) \geq 0$ *when* $(\hat{t}, \hat{n}) \in (0, T) \times \{0\}$, *and* $\bar{\Phi}(\hat{t}, \hat{n}) \geq h(\hat{n})$ *when* $(\hat{t}, \hat{n}) \in \{T\} \times [0, +\infty)$.

*(b) Viscosity sub-solution*

*Let $\bar{\Phi} \in USC(\bar{\Omega})$, which grows at most linearly with respect to the second argument. We refer to $\underline{\Phi}$ as a viscosity sub-solution of (34) if, for any $(\hat{t}, \hat{n}) \in \bar{\Omega}$ and $\varphi \in C^1(\bar{\Omega})$ such that $\underline{\Phi}(\hat{t}, \hat{n}) - \varphi(\hat{t}, \hat{n}) = \max_{(t,n) \in \bar{\Omega}} \{\underline{\Phi}(t,n) - \varphi(t,n)\}$, it follows that*

$$\frac{\partial \varphi}{\partial t}(\hat{t}, \hat{n}) - R(\hat{n})\hat{n}\frac{\partial \varphi}{\partial n}(\hat{t}, \hat{n}) + \tilde{H}\left(\hat{t}, \frac{\partial \varphi}{\partial n}(\hat{t}, \hat{n})\right) \geq 0 \quad \text{when} \quad (\hat{t}, \hat{n}) \in \Omega, \tag{37}$$

$\underline{\Phi}(\hat{t}, \hat{n}) \leq 0$ *when* $(\hat{t}, \hat{n}) \in (0, T) \times \{0\}$, *and* $\underline{\Phi}(\hat{t}, \hat{n}) \leq h(\hat{n})$ *when* $(\hat{t}, \hat{n}) \in \{T\} \times [0, +\infty)$.

*(c) Viscosity solution*

*A continuous function $\Phi : \bar{\Omega} \to \mathbb{R}$ is a viscosity solution if it is a viscosity super-solution in the sense of (a) and a viscosity sub-solution in the sense of (b).*

The viscosity solutions to the original HJBI equation are also defined.

## Definition 2

*The viscosity super-solutions, viscosity sub-solutions, and viscosity solutions to the HJBI equation (32) are defined by following **Definition 1**, where $\tilde{H}$ is formally replaced with $H$.*

At this point, we state the uniqueness result of the modified HJBI equation.

## Proposition 2

*The modified HJBI equation (34) admits at most one viscosity solution in the sense of **Definition 1**.*

We close this sub-section by analyzing the no-uncertainty limit $\mu \to +0$ of the HJBI equation. For this purpose, for $(t,z) \in [0,T] \times \mathbb{R}$, we set

$$\tilde{H}_0(t,z) = \frac{\alpha^2}{4(\gamma + \max\{0,z\})} \int_0^1 X_t(u) p(u) \mathrm{d}u. \tag{38}$$

A straightforward calculation shows that $\tilde{H}_\mu \to \tilde{H}_0$ locally uniformly in $(t,z) \in [0,T] \times \mathbb{R}$, which leads to the following stability result (e.g., Lemma 3.2 of [41]; Section 6 of [36]). In **Proposition 3**, the viscosity solutions to the HJBI equation in which $\tilde{H}$ is formally replaced with $\tilde{H}_0$ are defined analogously to those of **Definition 1**. This proposition states that the proposed HJBI equation that accounts for the model uncertainty is consistent with that without the uncertainty.

***Proposition 3***

*If $\{\Phi_\mu\}_{\mu>0}$ is a sequence of viscosity solutions in the sense of **Definition 1** for each fixed $\mu > 0$ that is bounded in each open set $\omega$ in $\Omega$. If the half-relaxed limit $\underline{\Phi}(t,n) = \liminf_{((t',n') \in \omega) \to (t,n),\, \mu \to +0} \Phi_\mu(t',n')$ (resp., $\overline{\Phi}(t,n) = \limsup_{((t',n') \in \omega) \to (t,n),\, \mu \to +0} \Phi_\mu(t',n')$) is continuous on $\overline{\Omega}$, it is a viscosity super-solution (resp., sub-solution) to the equation*

$$\frac{\partial \Phi}{\partial t} - R(n) n \frac{\partial \Phi}{\partial n} + \tilde{H}_0\left(t, \frac{\partial \Phi}{\partial n}\right) = 0, \tag{39}$$

*subject to the boundary and terminal conditions (13) and (14).*

### 3.3 Finite difference method
#### 3.3.1 Formulation

We present a monotone finite difference method that is fully explicit in time for the modified HJBI equation (34). The domain $\Omega$ is first discretized into a computational grid

$$\Omega_{\text{num}} = \{(i\Delta t, j\Delta n) | 0 \leq i \leq I, 0 \leq j \leq J\}, \tag{40}$$

where $\Delta t, \Delta n > 0$ are the grid sizes in time $t$ and population $n$, respectively. The degree of freedom of the grid is controlled by $I, J \in \mathbb{N}$. We assume that $I\Delta t = T$ as well as $J\Delta n = M$ with a constant $M > 0$. In practice, the parameter $M$ can be selected to be the maximum population at the initial time $t = 0$, which can be estimated from an ecological survey.

The discretized solution $\Phi$ to the HJBI equation (32) at the grid point $P_{i,j} := (i\Delta t, j\Delta n)$ is denoted as $\Phi_{i,j}$. Our finite difference method directly specifies the boundary condition as

$$\Phi_{i,0} = 0, \quad 0 \leq i \leq I-1 \tag{41}$$

and the terminal condition as

$$\Phi_{I,j} = h(j\Delta n), \quad 0 \leq j \leq J. \tag{42}$$

Subsequently, the HJBI equation itself is discretized at each $P_{i,j}$ ($0 \leq i \leq I-1$, $1 \leq j \leq J$) as follows:

$$\frac{\Phi_{i+1,j} - \Phi_{i,j}}{\Delta t} - R(j\Delta n) j\Delta n \frac{\Phi_{i+1,j} - \Phi_{i+1,j-1}}{\Delta n} + \tilde{H}\left(i\Delta t, \frac{\Phi_{i+1,j} - \Phi_{i+1,j-1}}{\Delta n}\right) = 0, \tag{43}$$

which can be rearranged as (recall that we are dealing with an equation to be solved backward in time)

$$\begin{aligned}\Phi_{i,j} &= \Phi_{i+1,j} - R(j\Delta n) j\Delta n \Delta t \frac{\Phi_{i+1,j} - \Phi_{i+1,j-1}}{\Delta n} + \tilde{H}\left(i\Delta t, \frac{\Phi_{i+1,j} - \Phi_{i+1,j-1}}{\Delta n}\right)\Delta t \\ &:= G(i, j, \Phi_{i+1,j}, \Phi_{i+1,j-1})\end{aligned} \tag{44}$$

In our numerical computation in **Section 4**, the integration with respect to $u$ that is involved in $\tilde{H}$ is carried out by a midpoint rule with uniformly distributed at 150 points ($u = (m-0.5)/150$, $m = 1, 2, 3, ..., 150$), which was preliminary found to be sufficient for our application. Using (41), (42), and (44), we can obtain the numerical solution $\Phi_{i,j}$ ($0 \leq i \leq I$, $0 \leq j \leq J$) from $i = I$ to $i = 0$ without resorting to solving any matrix inversion. However, it is necessary to take a sufficiently small $\Delta t$ for stability (See, next subsection).

***Remark 6*** Using (44), we do not require information outside the computational domain; that is, we do not need to specify any artificial boundary conditions along $j = J$. This is owing to the mathematical structure of the discretized and original HJBI equations that the information propagates from $n = 0$ to larger $n$ values.

### 3.3.2 Analysis of numerical method

We first prove the following stability result, and state that the non-negativity of the numerical solution is satisfied for a sufficiently small $\Delta t$. The proposed finite difference method can be stabilized, whereas the computational cost increases linearly with respect to $1/\Delta t$ owing to the fully explicit nature.

***Proposition 4***

*The non-negativity* $\Phi_{i,j} \geq 0$ *($0 \leq i \leq I$, $0 \leq j \leq J$) is satisfied if*

$$0 < \Delta t \leq \Delta n \left(R(M)M + \frac{\alpha^2}{4\gamma^2} K(1)\right)^{-1}. \tag{45}$$

Furthermore, we obtain an upper bound of the numerical solution that is uniform in time.

***Proposition 5***

*From (45), it follows that* $\Phi_{i,j} \leq h(M) + \dfrac{\alpha^2 T}{4\gamma} K(1)$ *(* $0 \leq i \leq I$, $0 \leq j \leq J$ *).*

**Propositions 4** and **5** demonstrate that the numerical solutions generated by the proposed finite difference method are uniformly stable for a sufficiently small time increment that satisfies (45). Furthermore, an increasing property of the numerical solutions with respect to $j$ is satisfied in the HJBI equations as a byproduct of the monotonicity. The proof of the proposition uses the argument of Lemma 5.2, 3(a) of [55]. This is a discrete counterpart of **Proposition 1**.

*Proposition 6*

*Under (45), the right-hand side of (44) increases with respect to both* $\Phi_{i+1,j}$ *and* $\Phi_{i+1,j-1}$. *Furthermore, it follows that* $\Phi_{i,j} \leq \Phi_{i,j+1}$ *(* $0 \leq i \leq I$, $0 \leq j \leq J-1$ *).*

*Remark 7* According to **Propositions 4–6**, scheme (44) is essentially the same as

$$\Phi_{i,j} = \Phi_{i+1,j} - R(j\Delta n) j\Delta n \Delta t \frac{\Phi_{i+1,j} - \Phi_{i+1,j-1}}{\Delta n} + H\left(i\Delta t, \frac{\Phi_{i+1,j} - \Phi_{i+1,j-1}}{\Delta n}\right)\Delta t, \qquad (46)$$

which is the naïve finite difference discretization of the original HJBI equation (32).

*Remark 8* One may expect that the numerical solutions that are generated by the finite difference method converge to a viscosity solution in the sense o**f Definition 1** locally uniformly in $(0,T) \times (0,M)$ if $\Delta t = \varsigma \Delta n$ is selected for a sufficiently small constant $\varsigma > 0$ that at least complies with (45). This follows from the classical convergence argument [38] because the scheme is monotone, stable, and clearly consistent, the Hamiltonian $H$ is Lipschitz continuous, and both the boundary and terminal conditions are Lipschitz continuous as well as bounded. See also Theorem 3.15 and Section 4.1 of [56]. Therefore, the numerical solutions yield the existence of a viscosity solution; hence, they converge to the unique viscosity solution by **Proposition 2**. More specifically, we prepare the monotone bilinear interpolation $\Phi_\Delta$ of the numerical solutions: for each $(t,n) \in S_{i,j} := [i\Delta t, (i+1)\Delta t] \times [j\Delta n, (j+1)\Delta n]$ ( $0 \leq i \leq I-1, 0 \leq j \leq J-1$ ), we set

$$\Phi_\Delta(t,n) := \frac{(i+1)\Delta t - t}{\Delta t} \frac{(j+1)\Delta n - n}{\Delta n} \Phi_{i,j} + \frac{t - i\Delta t}{\Delta t} \frac{(j+1)\Delta n - n}{\Delta n} \Phi_{i+1,j}$$
$$+ \frac{(i+1)\Delta t - t}{\Delta t} \frac{n - j\Delta n}{\Delta n} \Phi_{i,j+1} + \frac{t - i\Delta t}{\Delta t} \frac{n - j\Delta n}{\Delta n} \Phi_{i+1,j+1}. \qquad (47)$$

By construction, this $\Phi_\Delta$ is continuous over $[0,T] \times [0,M]$, and satisfies the bounds in **Propositions 4** and **5**. According to the Ascoli–Arzelà theorem and **Propositions 4** and **5**, along with the selection of the discretization parameters ( $\Delta t = \varsigma \Delta n$ ), there exists a subsequence of $\Phi_\Delta$ that converges locally uniformly in $(0,T)\times(0,M)$.

***Remark 9*** Relating to **Remark 8**, we find that the convergence that preserves the increasing property in $n$ appears to be far trivial. Such a proof can be found in the proof of Theorem 5.5 of [55]. However, their proof fails in our case owing to the lack of no increasing or decreasing property with respect to $t$. Indeed, the numerical solutions that are obtained in **Section 4** suggest that the solution to the (modified) HJBI equation is non-monotone with respect to $t$.

## 4. Application

### 4.1 Study site

For the model application, we used growth data of *P. altivelis* that were collected in the Hii River, which is one of the largest rivers in the San-in area, Shimane Prefecture, Japan, with the help of the Hii River Fishery Cooperative (HRFC) which authorizes inland fishery resources in the mid-to-up-stream reaches of this river. The authors had been communicating with the HRFC since 2015, and previously investigated the environment and fisheries of the Hii River (e.g., see [57,58] and the references therein). The fish *P. altivelis* is a major inland fishery resource that contributes to the regional environment, ecosystem, and culture and recreation. Tomozuri (fishing with decoys [59]) and Toami (casting a net [60]) are major fishing methods for catching the fish.

*P. altivelis* has a one-year life history such that they spawn in a river during autumn, and the spawn fishes flow down along the river to the downstream sea or a reservoir, overwinter in the sea or a reservoir, and migrate towards the river in the spring to mature during the summer (e.g., [61]). The harvesting period of *P. altivelis* in the Hii River usually starts on July 1 ($t = 61$ (days) below) and ends in October to November ($T = 181$ (days) below), where May 1 was selected as the reference point of $t = 0$ (days) following previous studies (e.g., [57]), which corresponds to the growing season of the fish. Therefore, the computational time domain is $(61,181)$ (days) with the length of 120 (days). The population dynamics of *P. altivelis* in the Hii River has not been clarified to date. In particular, no estimate has been made for the overall river system. However, it has recently been estimated that $O(10^4)$ juveniles of *P. altivelis* were found in the upstream main branch of the Sakura–Orochi reservoir, which creates the largest dam, namely the Obara Dam, in the Hii River [62]. Thus, we focused on the application of the proposed model to the river based on the assumption that the maximum population at the beginning of July was $O(10^4)$.

### 4.2 Computational conditions

Sampling surveys of the body weights of *P. altivelis* individuals are carried out by a union member of the HRFC from Summer to Autumn each year, which we used to estimate the model parameters of $X$. We used the most recent data that were collected during 2021 and 2022, as illustrated in **Figure 1**. Note that different individuals were caught at different observation times.

The fitting method of the uncertain logistic model (5) is based on the maximum likelihood method of Dorini et al. [33], whose conclusion demonstrated that the uniform distribution should be assumed with a minimum data requirement provided that the upper bound $K(1)$ and lower bound $K(0)$ are given. Within this framework, we assumed that only $K$ is uncertain, leading to the parameters that were estimated in the uncertain logistic model being $x$, $r$, and $\bar{K}, \underline{K}$. For each year, we attempted to find the best values of $(x, r, \bar{K}, \underline{K})$ to minimize the common least-squares error loss function between the observed model ($X_{\text{obs},i}$) and fitted model ($X_{\text{theor},i}$):

$$\frac{1}{N_{\text{obs}}} \sum_{i=1}^{N_{\text{obs}}} (X_{\text{obs},i} - X_{\text{theor},i})^2, \tag{48}$$

where $N_{\text{obs}}$ is the total number of observations (40 in 2021 and 59 in 2022), $X_{\text{obs},i}$ is the observed weight at the $i$ th observation time, and $X_{\text{theor},i}$ is the theoretical mean body weight at the $i$ th observation time. The minimizer of the loss function (48) was determined using a trial-and-error approach such that we selected the best quadruple $(x, r, \bar{K}, \underline{K})$ to minimize (48) among the parameter values of $r = 0.020\text{-}0.050$ (1/day): increments of 0.001 (1/day), $x = 5\text{-}15$ (g): increments of 1 (g), $\underline{K} = 1\text{-}301$ (g): increments of 1 (g), and $\bar{K} = 1\text{-}301$ (g): increments of 1 (g) with the constraint $\bar{K} \geq \underline{K} + 1$ (g).

**Figure 2** compares the collected data and the fitted model of the body weights of individuals for each year. **Table 1** summarizes the fitted parameter values for each year. The fitted results tracked the observed data reasonably. The parameter values in **Table 1** suggest that between 2021 and 2022, the fitted model of *P. altivelis* in 2021 had larger individual differences than that of 2022, with a higher growth speed. Furthermore, the individuals were estimated to be smaller than those in 2022 near May 1; that is, around the season of upstream migration. The reason for these differences can be attributed to different hydrological and meteorological conditions as well as the associated ecosystem dynamics, which are difficult to clarify based only on the available data and are beyond the scope of this paper. Nevertheless, the fitted results imply that the uncertain logistic model could effectively capture the growth dynamics of the fish for different years.

The other parameter values were similar for the two years. The population dynamics (1) was normalized to $10^4$ individuals considering the discussion in the previous sub-section. Therefore, in the computation, the individuals were counted in units of $10^4$. Under this normalization, the maximum domain size $M$ in the $n$ direction was set as 10, corresponding to $10^5$ individuals. The mortality $R$ was set to 0.01 (1/day) for simplicity considering previous identification results [63]. The parameters in the objective function were set to $\alpha = 0.05$ (1/day), $\gamma = 0.1$ (1/day), $\mu = 0.01$ (1/day), and $h(n) = 0$ (no sustainability concern) or $h(n) = 1.5\min\{n, M\}$ (with the sustainability concern); these values had been preliminarily explored so that realistic sample paths of the controlled resource dynamics could be obtained as well as we can visually demonstrate their influences on the optimal harvesting policy as well as the worst-case uncertainty. Particularly, we have heuristically chosen the parameter value of $\mu$ so that the uncertainty affects the computational results while it will not lead to a too small estimate of the body weight (e.g., see **Figure 7**). This parameter is user-specific and quantifies how much the decision-maker concerns the uncertainty. Its estimation will need another study about decision-making mechanisms inside humans, which is beyond the scope of this study. The computational resolutions were $I = 120\ 000$ and $J = 1\ 000$, which were found to be sufficient for the computation below.

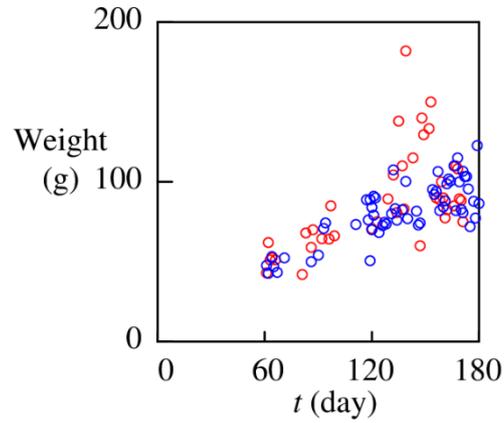

**Figure 1.** Collected data of body weights of individuals in 2021 (red) and 2022 (blue).

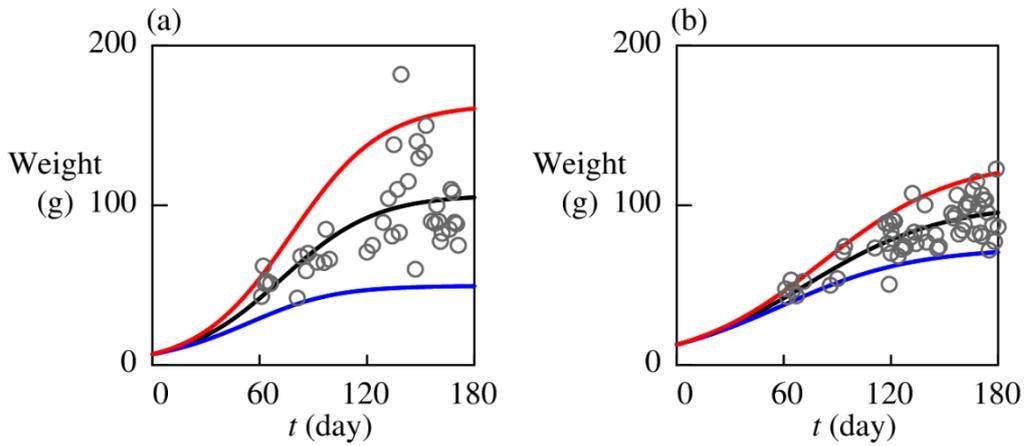

**Figure 2.** Comparison of collected data and fitted model of body weights of individuals in **(a)** 2021 and **(b)** 2022. In each plot, the circles denote observations, the black curve indicates the theoretical mean, the red curve denotes the theoretical mean + standard deviation, and the blue curve indicates the theoretical mean − standard deviation.

**Table 1.** Fitted parameter values of the uncertain logistic model in each year.

|  | 2021 | 2022 |
|---|---|---|
| $x$ (g) | 6.8 | 12.8 |
| $r$ (g/day) | 0.040 | 0.027 |
| $\overline{K}$ (g) | 205 | 149 |
| $\underline{K}$ (g) | 8 | 53 |

## 4.3 Computational results

**Figures 3(a)** and **3(b)** depict the computed value functions $\Phi$ for the years 2021 and 2022, respectively, without the sustainability concern. Similarly, **Figures 4(a)** and **4(b)** show the computed value functions $\Phi$ for the years 2021 and 2022, respectively, with the sustainability concern. The computed value functions appeared to be continuous and smooth inside the domain in both years, regardless of the existence of the sustainability concern, thereby supporting the continuity assumption of the mathematical analysis. The computed $\Phi$ in 2022 was larger than that in 2021 in most parts of the computational domain, suggesting that the smaller individual differences led to a larger net utility under the assumed condition.

**Figures 5** and **6** depict the optimal harvesting policy and several controlled paths of the populations corresponding to **Figures 3** and **4**, respectively. The controlled paths were computed by numerically integrating the population dynamics (1) backwards in time based on the optimal harvesting policy that was obtained via the HJBI equation. This approach is justified in our case because the population dynamics does not contain stochastic noise. The set of terminal values of the backward paths is the same in each panel of **Figures 5** and **6**. A comparison of **Figures 5** and **6** reveals that accounting for the sustainability of the fish population critically affected the controlled population dynamics. According to **Figure 5**, for each year, the backward computed paths indicate that, given a terminal value, the population dynamics to achieve a positive terminal value required initial populations that were significantly larger than $10^5$ individuals. In contrast, **Figure 6** suggests that populations of the order of $O(10^4)$ archived the positive terminal values with the same order. Indeed, the optimal harvesting rates exhibited one order of difference between the cases with and without the sustainability concern. Although the required total number of populations to sustain the life cycles of the fish is dependent on the specific relationship between the body size and total number of eggs, as discussed for other fish species (e.g., [64,65]), the proposed mathematical model suggests that considering the terminal utility is key to the sustainable fisheries management of *P. altivelis* in the Hii River.

We also analyzed the growth dynamics under the uncertainty, which is a key element in our mathematical model. **Figure 7** presents paths of the distorted mean body weights (6) that were subject to different levels of uncertainty aversion for 2021 and 2022 with the sustainability concern. Note that if a classical deterministic logistic model is used, we formally only obtain $X = \bar{X}$, and hence, there will be no uncertainty effects in the control problem. In this study, the distorted $\bar{X}$ for each year was almost the same among the various sample paths with different terminal conditions, and the difference was at most several grams, so we only present the results of the growth dynamics corresponding to the terminal population of $10^4$ individuals (i.e., $n=1$ in the computation). The effects of the uncertainty aversion in terms of the weighting coefficient $\mu$ between the two years is clear, which demonstrates that critically different mean growth curves were generated for the same value of $\mu$, particularly for the smaller $\mu$ representing the larger uncertainty aversion. The larger potential individual differences in the model for 2021 were more susceptible to the model uncertainty, and an overly pessimistic prediction yielded an unrealistically small mean body weight (red curve in **Figure 7(a)**). One means of preventing model uncertainty to avoid such a pessimistic result is to collect the information at the beginning of the growth period, particularly before the harvesting season, which is the time 0 to 60 in our setting. Our model suggests that although such a survey would require high monetary and labor costs because the growth conditions are different among different years as demonstrated in this paper, its outcome can persist during the harvesting period to provide a less pessimistic prediction of the mean growth of the individual fishes. The computational also results suggest that the mean growth under the uncertainty approached that without the uncertainty, which is consistent with the theoretical results for the value function in **Proposition 3**.

We finally analyze the case with a non-uniform $p$ as our mathematical framework theoretically covers such a case. We focus on the beta distribution with parameters $a,b > 1$; that is, $p(u) = p_0 u^{a-1}(1-u)^{b-1}$ with a normalization constant $p_0$. The beta distribution is right-skewed, non-skewed, and left-skewed for $a < b$, $a = b$, and $a > b$, respectively (**Figure 8**). The data in the year 2021 is used here because of the larger uncertainty than 2022, with which the influences if the skewness can be more clearly visible. The right-skewed case seems to be the most qualitatively reasonable because the observed body weight data of *P. altivelis* in the past years were right-skewed (e.g., [11,48,57]). We consider all the three cases by specifying $(a,b) = (2,5)$, $(2,2)$, and $(5,2)$. The parameter values other than those in $p$ are the same with the computational cases presented above, and we set the uncertainty aversion as $\mu = 0.01$ (1/day). **Figure 9** shows the computed optimal harvesting policy and several controlled paths of populations in 2021 with sustainability concern for the right-skewed, non-skewed, and left-skewed cases. The skewness of the body weight distribution affects the optimal harvesting policy as well as the controlled population dynamics so that the smaller harvesting rate is optimal for a more right-skewed case. This is due to the modelling assumption that harvesting a larger individual is more profitable, meaning that too fast resource exploitation is less efficient for such a case. **Figure 10** presents computed paths of distorted mean body weights with sustainability concern, demonstrating that, even under the uncertainty, the left-skewed case predicts the larger mean body weight than that without uncertainty that is due to the distributional shape of $p$ concentrating more on the right-half of $0 \leq u \leq 1$ as shown in **Figure 8**. Accordingly, the non-skewed and right-skewed cases yield more pessimistic results than the left-skewed case. A comparison between **Figure 7(a)** and **Figure 10** suggests that using the beta distribution of the non-skewed case yields a more pessimistic result than the uniform case although the beta and uniform distributions share the same mean. The computational results thus suggest that the mean growth with the beta distribution, provided it is non- or right-skewed, is more pessimistic than the uniform distribution. Although the investigations of more realistic $p$ remains a future topic, its influences were qualitatively suggested by our numerical analysis.

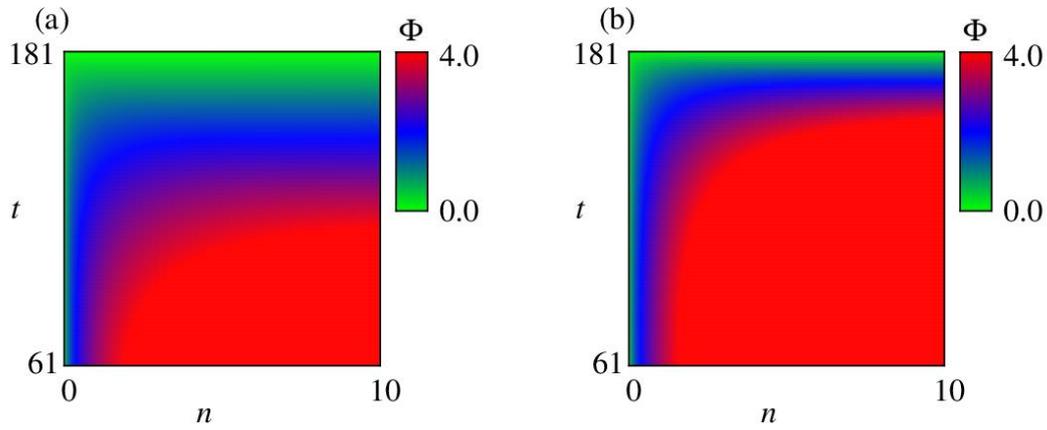

**Figure 3.** Computed value functions for **(a)** 2021 and **(b)** 2022 without sustainability concern.

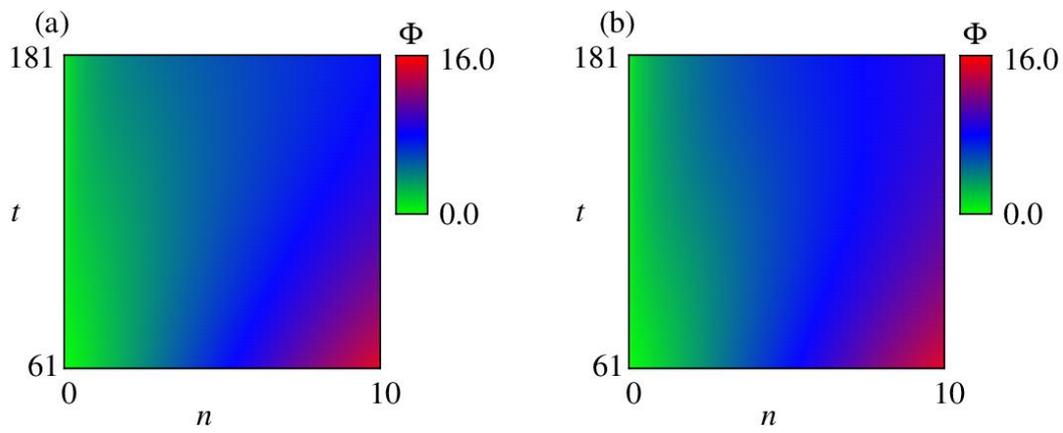

**Figure 4.** Computed value functions for **(a)** 2021 and **(b)** 2022 with sustainability concern.

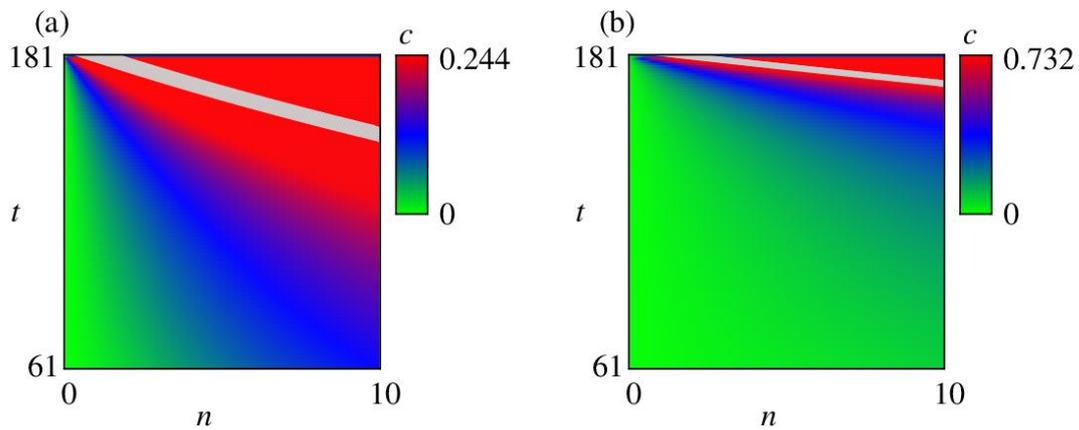

**Figure 5.** Computed optimal harvesting policy and several controlled paths of populations for **(a)** 2021 and **(b)** 2022 without sustainability concern. All computed controlled paths did not cross but were very close to one another.

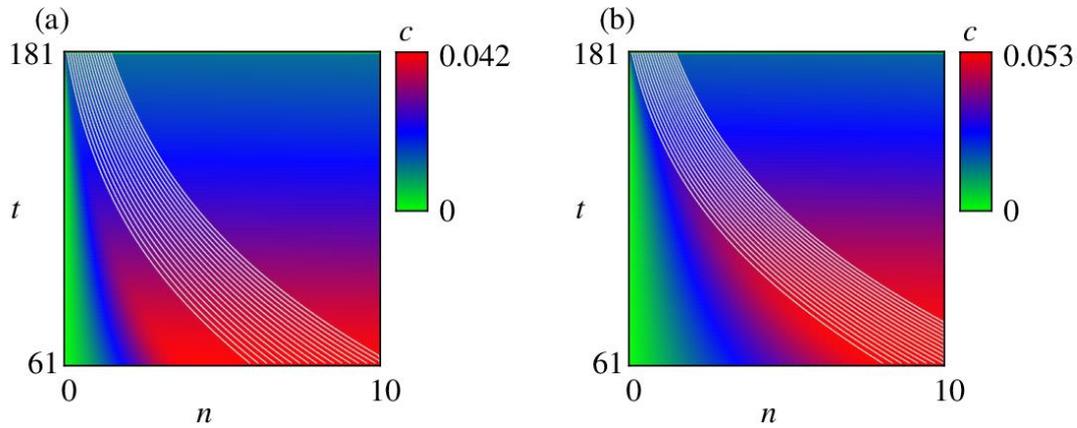

**Figure 6.** Computed optimal harvesting policy and several controlled paths of populations for **(a)** 2021 and **(b)** 2022 with sustainability concern. All computed controlled paths did not cross.

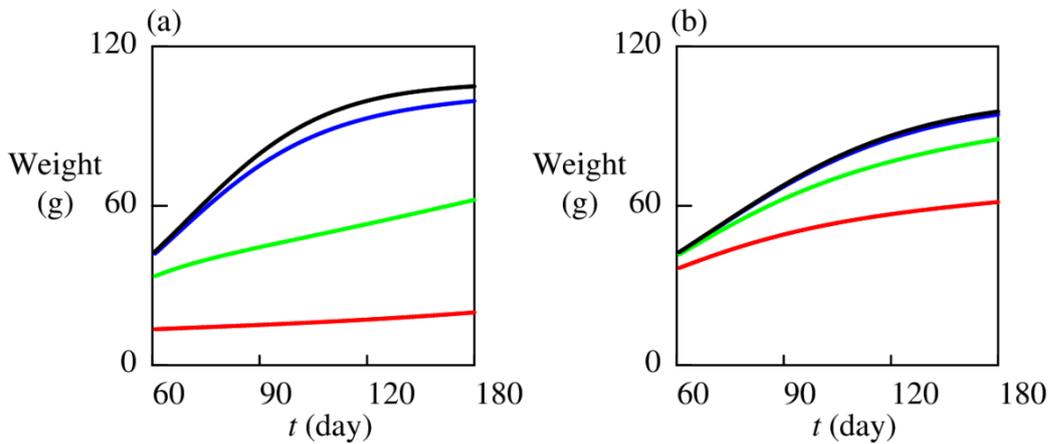

**Figure 7.** Computed paths of distorted mean body weights subject to different levels of uncertainty aversion for 2021 **(a)** without sustainability concern and **(b)** with sustainability concern. The values of the uncertainty aversion $\mu$ were $+\infty$ (black, no uncertainty), 0.1 (blue), 0.01 (green), and 0.001 (red).

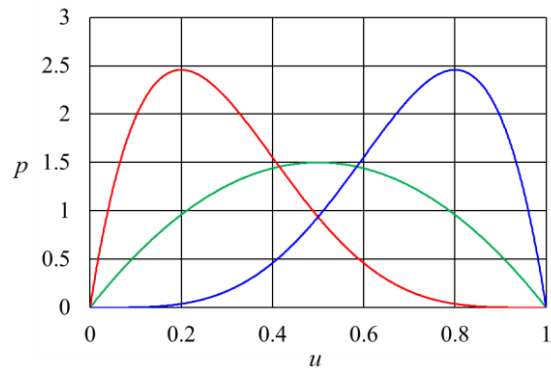

**Figure 8.** The beta distributions for the right-skewed case $(a,b)=(2,5)$ (red), non-skewed case $(a,b)=(2,2)$ (green), and left-skewed case $(a,b)=(5,2)$ (blue).

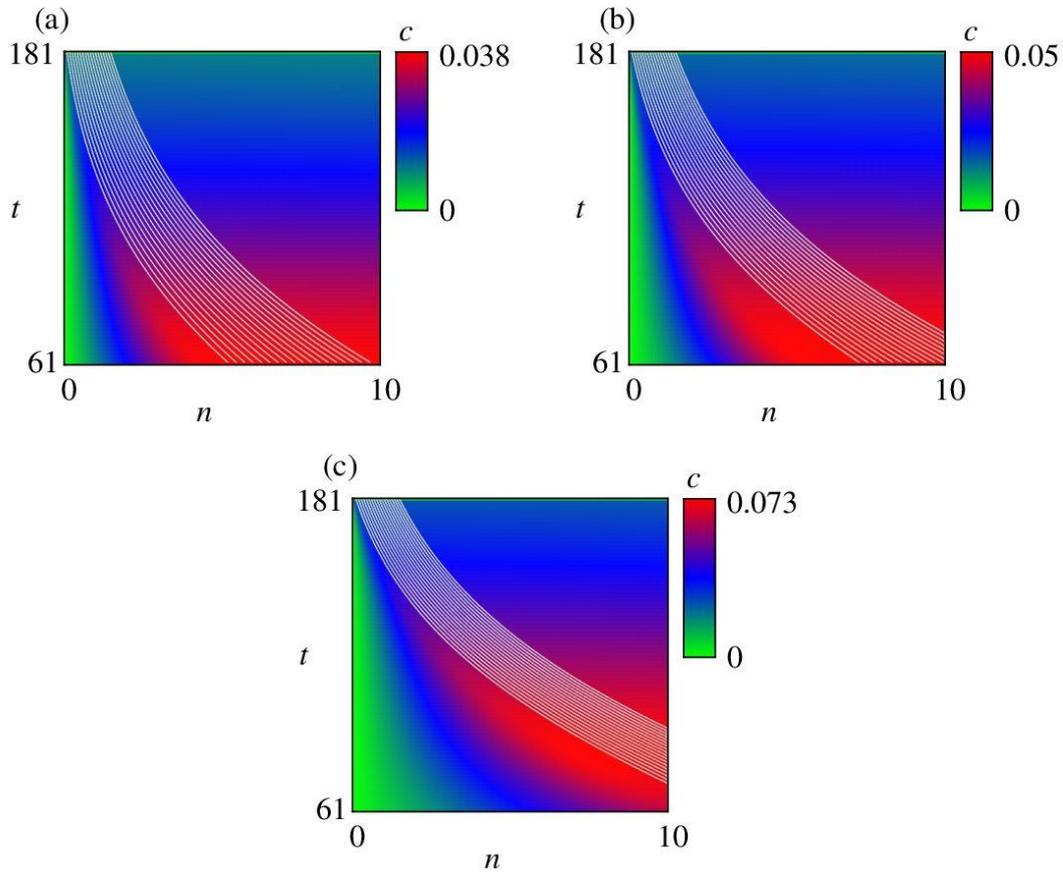

**Figure 9.** Computed optimal harvesting policy and several controlled paths of populations in 2021 with sustainability concern for the right-skewed case (a) $(a,b)=(2,5)$, non-skewed case (b) $(a,b)=(2,2)$, and left-skewed case (c) $(a,b)=(5,2)$. All computed controlled paths did not cross.

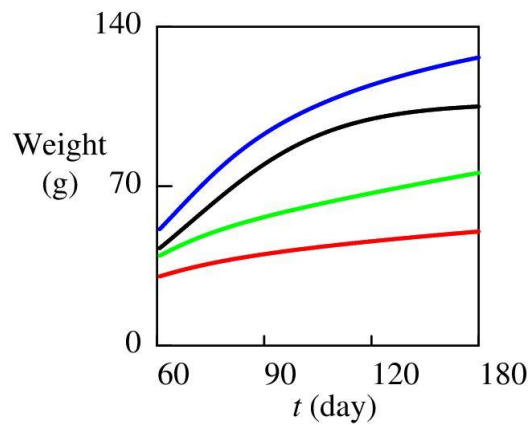

**Figure 10.** Computed paths of distorted mean body weights with sustainability concern for the right-skewed case $(a,b)=(2,5)$ (red), non-skewed case $(a,b)=(2,2)$ (green), and left-skewed case $(a,b)=(5,2)$ (blue).

## 5. Conclusions

This study proposed an optimal control problem for harvesting biological resources with physiological heterogeneity under model uncertainty. The Radon–Nikodym derivative was effectively employed to evaluate the uncertainty of the heterogeneity, and the resulting HJBI equation was successfully obtained in a tractable form. The mathematical analysis results characterized the shape and viscosity properties of the HJBI equation. A fully explicit and monotone finite difference method was presented and implemented in a real problem using the collected data. The optimal harvesting policy, which is dependent on the control objective and degree of uncertainty aversion, was successfully computed.

The proposed model can be applied to any other biological resources having physiological heterogeneity, such as a fish species other than *P. altivelis* [66]. We consider that each fish species can be characterized by the growth curve and population decrease/increase because they are used to estimate the biomass. Once the dynamics of the two quantities of the target fish have been identified, the proposed control setting can be applied to address its resource management problem. The HJBI equation for the management problem will then be obtained. Adding a noise process such as Brownian motions to the population dynamics is theoretically possible, while the difficulty lies in the estimation of the volatility (i.e., noise intensity) from field data as it would need time series of growth paths of individuals, which are difficult to track in practice. Nevertheless, one may be able to estimate the volatility based on the assumption that the growth rate $r$ should be larger than growth declines due to stochastic fluctuations in some sense. For example, one may assume the increasing property of the second moment of the body weight (e.g., fishes will continue to grow in some sense). In our context, if we assume the no-uncertainty case where the body weight is a time-dependent quantity $X_t$ ($K, r$ are constants) for simplicity and also assume the stochastic logistic model with the fluctuation term of the form $\varpi X_t \mathrm{d}B_t$ with the volatility $\varpi > 0$ and the 1-D standard Brownian motion $B_t$ (e.g., [48]), then the assumption stated above suggests that $\varpi$ should satisfy $r > \dfrac{\varpi^2}{2}$ by linearizing the nonlinearity of the logistic growth. This gives a rough upper bound of $\varpi$.

Several topics will be addressed in our future works. Resource harvesting problems in an infinite horizon [67] and that under an economic equilibrium [68] can be investigated based on the proposed approach by setting the terminal time to be sufficiently large. The size-dependent mortality [69] was not considered in this study owing to the increased complexity of the HJBI equation, such that the built-in optimization problem, namely the inf part, possibly becomes intractable and/or ill posed; this will be studied in the future. A more sophisticated numerical method, such as the fitted finite volume method [70], will be necessary to consider the noise-driven resource dynamics. The extension of the proposed mathematical approach to a multi-site problem [71] will also be an interesting future direction, wherein the efficient numerical computation of the associated optimality equation could be a potential challenge. Population dynamics with a delay sometimes arises in applications, and this can also be handled using dynamic programming, although the problem dimension becomes infinite [72] in such cases. Agent-based models as alternative representations of the heterogeneity are also interesting subjects that should consider high-dimensional dynamical systems. For example, it will be interesting to extend the proposed model so that the social interactions [20, 21] can be accounted for. A model reduction will be necessary to analyze such advanced problems to make the problem computationally feasible. Incorporating a statistical filter [48] to improve the uncertainty estimate will be another option to extend the proposed model. We focused on a problem in a natural environment, while the proposed approach also applies to that in an artificial system such as an aquacultural system [73].

**Appendix A: Proofs of propositions and lemmas**

*Proof of Lemma 1*

The proof is owing to the inequality

$$\frac{\alpha^2}{4(\gamma+z)^2}\int_0^1 X_t(u)\phi(u)p(u)\mathrm{d}u \leq \frac{\alpha^2}{4(\gamma+z)^2}\int_0^1 K(1)\phi(u)p(u)\mathrm{d}u$$

$$= \frac{\alpha^2}{4\gamma^2}K(1)\int_0^1 \phi(u)p(u)\mathrm{d}u \qquad (49)$$

$$= \frac{\alpha^2}{4\gamma^2}K(1)$$

and the equality due to **Assumption 1**:

$$\hat{c}(t,z,\phi) = \max\left\{0, \min\left\{\frac{\alpha^2}{4(\gamma+z)^2}\left(\int_0^1 X_t(u)\phi(u)p(u)\mathrm{d}u\right), \bar{c}\right\}\right\}$$

$$= \min\left\{\frac{\alpha^2}{4(\gamma+z)^2}\left(\int_0^1 X_t(u)\phi(u)p(u)\mathrm{d}u\right), \bar{c}\right\} \qquad (50)$$

$$= \frac{\alpha^2}{4(\gamma+z)^2}\int_0^1 X_t(u)\phi(u)p(u)\mathrm{d}u$$

□

*Proof of Proposition 1*

We first present a comparison result of the controlled system:

$$N_s^{(t,n_2,c,\phi)} \geq N_s^{(t,n_1,c,\phi)}, \quad t \leq s \leq \min\{\tau_{t,n_1}, \tau_{t,n_2}\}. \qquad (51)$$

If $n_1 = n_2$ or $n_1 = 0$, (51) is trivial. Therefore, we assume that $0 < n_1 < n_2$. We fix $c \in \mathbb{A}$ and $\phi \in \mathbb{B}$, and obtain

$$N_s^{(t,n,c,\phi)} = N_t^{(t,n,c,\phi)} \exp\left(-\int_t^s R\left(N_u^{(t,n,c,\phi)}\right)\mathrm{d}u\right) - \int_t^s \exp\left(-\int_v^s R\left(N_v^{(t,n,c,\phi)}\right)\mathrm{d}v\right)c_u\mathrm{d}u, \quad 0 < s < \tau_{t,n} \qquad (52)$$

for $n = n_1, n_2$. We set $\tau' = \inf\{s \geq t | N_s^{(t,n_2,c,\phi)} = N_s^{(t,n_1,c,\phi)}\}$. Owing to the continuity of (52) and $0 < n_1 < n_2$, this $\tau'$ must be larger than $t$. Furthermore, it must be smaller than $\max\{\tau_{t,n_1}, \tau_{t,n_2}\}$. We also obtain $\min\{\tau_{t,n_1}, \tau_{t,n_2}\} > 0$.

For $t < s < \min\{\tau', \min\{\tau_{t,n_1}, \tau_{t,n_2}\}\}$, we obtain

$$\frac{\mathrm{d}N_s^{(t,n_i,c,\phi)}}{\mathrm{d}s} = -R\left(N_s^{(t,n_i,c,\phi)}\right)N_s^{(t,n_i,c,\phi)} - c_s, \quad i = 1,2, \qquad (53)$$

and hence, due to $N_s^{(t,n_2,c,\phi)} \geq N_s^{(t,n_1,c,\phi)}$, the estimate

$$\begin{aligned}
\frac{d\left(N_s^{(t,n_2,c,\phi)} - N_s^{(t,n_1,c,\phi)}\right)}{ds} &= -R\left(N_s^{(t,n_2,c,\phi)}\right)N_s^{(t,n_2,c,\phi)} - c_s - \left(-R\left(N_s^{(t,n_1,c,\phi)}\right)N_s^{(t,n_1,c,\phi)} - c_s\right) \\
&= -R\left(N_s^{(t,n_2,c,\phi)}\right)N_s^{(t,n_2,c,\phi)} + R\left(N_s^{(t,n_1,c,\phi)}\right)N_s^{(t,n_1,c,\phi)} \\
&\geq -L'\left|N_s^{(t,n_2,c,\phi)} - N_s^{(t,n_1,c,\phi)}\right| \\
&= -L'\left(N_s^{(t,n_2,c,\phi)} - N_s^{(t,n_1,c,\phi)}\right)
\end{aligned} \qquad (54)$$

with a positive constant $L'$ that is independent from $s, n_1, n_2, c, \phi$ owing to the boundedness and Lipschitz continuity of $R$. Subsequently, an application of the classical Gronwall inequality demonstrates that

$$\begin{aligned}
N_s^{(t,n_2,c,\phi)} - N_s^{(t,n_1,c,\phi)} &\geq \left(N_t^{(t,n_2,c,\phi)} - N_t^{(t,n_1,c,\phi)}\right)\exp(-L'(s-t)) \\
&= (n_2 - n_1)\exp(-L'(s-t))
\end{aligned} \qquad (55)$$

If $\tau' < \min\{\tau_{t,n_1}, \tau_{t,n_2}\}$, we obtain

$$0 = N_{\tau'}^{(t,n_2,c,\phi)} - N_{\tau'}^{(t,n_1,c,\phi)} \geq (n_2 - n_1)\exp(-L'(\tau'-t)); \qquad (56)$$

that is, $0 \geq \exp(-L'(\tau'-t))$ (recall that $n_2 > n_1$), which is impossible unless $\tau' = +\infty$: a contradiction. Therefore, we must have (51). As a byproduct, the continuity of (52) indicates that we must obtain $\tau_{t,n_1} \leq \tau_{t,n_2}$. Consequently, we can extend (51) as follows:

$$N_s^{(t,n_2,c,\phi)} \geq N_s^{(t,n_1,c,\phi)}, \quad t \leq s \leq T. \qquad (57)$$

At this point, we fix $\phi \in \mathbb{B}$ and $c \in \mathbb{A}$ such that $c_s = 0$ for $s \geq \tau_{t,n_1}$. In particular, we have $c_s = 0$ for $\min\{T, \tau_{t,n_1}\} < s < \min\{T, \tau_{t,n_2}\}$. Such a control is admissible for both $n = n_2, n_1$. Subsequently, we obtain

$$J(t, n_i, c, \phi) = \int_t^T \left(\alpha(c_s \bar{X}_{s,\phi})^\beta - \gamma c_s\right)ds + h\left(N_T^{(t,n_i,c,\phi)}\right) + \mu\int_t^T \mathbb{D}(\phi_s)ds \qquad (58)$$

and hence,

$$\begin{aligned}
J(t,n_2,c,\phi) - J(t,n_1,c,\phi) &= \int_t^T\left(\alpha(c_s\bar{X}_{s,\phi})^\beta - \gamma c_s\right)ds + h\left(N_T^{(t,n_2,c,\phi)}\right) + \mu\int_t^T\mathbb{D}(\phi_s)ds \\
&\quad - \left(\int_t^T\left(\alpha(c_s\bar{X}_{s,\phi})^\beta - \gamma c_s\right)ds + h\left(N_T^{(t,n_1,c,\phi)}\right) + \mu\int_t^T\mathbb{D}(\phi_s)ds\right) \\
&= \int_t^T\left(\alpha(c_s\bar{X}_{s,\phi})^\beta - \gamma c_s\right)ds - \int_t^T\left(\alpha(c_s\bar{X}_{s,\phi})^\beta - \gamma c_s\right)ds \\
&\quad + h\left(N_T^{(t,n_2,c,\phi)}\right) - h\left(N_T^{(t,n_1,c,\phi)}\right) \\
&= \int_t^{\min\{T,\tau_{t,n_2}\}}\left(\alpha(c_s\bar{X}_{s,\phi})^\beta - \gamma c_s\right)ds - \int_t^{\min\{T,\tau_{t,n_1}\}}\left(\alpha(c_s\bar{X}_{s,\phi})^\beta - \gamma c_s\right)ds \quad (59)\\
&\quad + h\left(N_T^{(t,n_2,c,\phi)}\right) - h\left(N_T^{(t,n_1,c,\phi)}\right) \\
&= \int_{\min\{T,\tau_{t,n_1}\}}^{\min\{T,\tau_{t,n_2}\}}\left(\alpha(c_s\bar{X}_{s,\phi})^\beta - \gamma c_s\right)ds + h\left(N_T^{(t,n_2,c,\phi)}\right) - h\left(N_T^{(t,n_1,c,\phi)}\right) \\
&= h\left(N_T^{(t,n_2,c,\phi)}\right) - h\left(N_T^{(t,n_1,c,\phi)}\right) \\
&\geq 0
\end{aligned}$$

owing to the increasing property of $h$. Consequently, for each $\phi \in \mathbb{B}$, and fixing $c \in \mathbb{A}$ such that $c_s = 0$ for $s \geq \tau_{t,n_1}$, we obtain

$$J(t,n_2,c,\phi) \geq J(t,n_1,c,\phi). \tag{60}$$

The left-hand side of (60) is not larger than $\sup_{c \in \mathbb{A}} J(t,n_2,c,\phi)$, and hence

$$\sup_{c \in \mathbb{A}} J(t,n_2,c,\phi) \geq J(t,n_1,c,\phi). \tag{61}$$

Then, we obtain (the restriction $c_s = 0$ for $s \geq \tau_{t,n_1}$ does not affect taking the supremum of the right-hand side below)

$$\sup_{c \in \mathbb{A}} J(t,n_2,c,\phi) \geq \sup_{c \in \mathbb{A}} J(t,n_1,c,\phi). \tag{62}$$

In the same manner, we have

$$\sup_{c \in \mathbb{A}} J(t,n_2,c,\phi) \geq \inf_{\phi \in \mathbb{B}} \sup_{c \in \mathbb{A}} J(t,n_1,c,\phi), \tag{63}$$

and hence,

$$\inf_{\phi \in \mathbb{B}} \sup_{c \in \mathbb{A}} J(t,n_2,c,\phi) \geq \inf_{\phi \in \mathbb{B}} \sup_{c \in \mathbb{A}} J(t,n_1,c,\phi), \tag{64}$$

which is the desired inequality (31).

□

*Proof of Lemma 2*

The uniform continuity of $H$ follows from the uniform continuity of $e^{-\frac{\alpha^2}{4\mu(\gamma+z)}X_t(u)}$, and hence, that of $\int_0^1 e^{-\frac{\alpha^2}{4\mu(\gamma+z)}X_t(u)} p(u)\mathrm{d}u$, which is further combined with the lower bound

$$\begin{aligned}
-\mu \ln\left(\int_0^1 e^{-\frac{\alpha^2}{4\mu(\gamma+z)}X_t(u)} p(u)\mathrm{d}u\right) &\geq -\mu \ln\left(\int_0^1 e^{-\frac{\alpha^2}{4\mu\gamma}K(0)} p(u)\mathrm{d}u\right) \\
&= -\mu \ln\left(e^{-\frac{\alpha^2}{4\mu\gamma}K(0)}\right) \\
&= \frac{\alpha^2}{4\gamma}K(0) \\
&> 0
\end{aligned} \tag{65}$$

For the Lipschitz continuity, we obtain

$$\frac{\partial H}{\partial z}(t,z) = \frac{\partial}{\partial z}\left\{-\mu \ln\left(\int_0^1 e^{-\frac{\alpha^2}{4\mu(\gamma+z)}X_t(u)} p(u)\,du\right)\right\}$$

$$= -\mu \frac{1}{\int_0^1 e^{-\frac{\alpha^2}{4\mu(\gamma+z)}X_t(u)} p(u)\,du} \frac{\partial}{\partial z}\int_0^1 e^{-\frac{\alpha^2}{4\mu(\gamma+z)}X_t(u)} p(u)\,du$$

$$= -\mu \frac{1}{\int_0^1 e^{-\frac{\alpha^2}{4\mu(\gamma+z)}X_t(u)} p(u)\,du} \int_0^1 \frac{\partial}{\partial z} e^{-\frac{\alpha^2}{4\mu(\gamma+z)}X_t(u)} p(u)\,du \quad , \quad (66)$$

$$= -\mu \frac{1}{\int_0^1 e^{-\frac{\alpha^2}{4\mu(\gamma+z)}X_t(u)} p(u)\,du} \int_0^1 \frac{\alpha^2}{4\mu(\gamma+z)^2} X_t(u) e^{-\frac{\alpha^2}{4\mu(\gamma+z)}X_t(u)} p(u)\,du$$

$$= -\frac{\alpha^2}{4(\gamma+z)^2} \frac{1}{\int_0^1 e^{-\frac{\alpha^2}{4\mu(\gamma+z)}X_t(u)} p(u)\,du} \int_0^1 X_t(u) e^{-\frac{\alpha^2}{4\mu(\gamma+z)}X_t(u)} p(u)\,du$$

where the order of the partial differentiation and integration can be exchanged owing to the smoothness and uniform boundedness of $e^{-\frac{\alpha^2}{4\mu(\gamma+z)}X_t(u)}$ for $z \geq 0$. We then proceed as follows:

$$\left|\frac{\partial H}{\partial z}(t,z)\right| = \frac{\alpha^2}{4(\gamma+z)^2} \frac{1}{\int_0^1 e^{-\frac{\alpha^2}{4\mu(\gamma+z)}X_t(u)} p(u)\,du} \int_0^1 X_t(u) e^{-\frac{\alpha^2}{4\mu(\gamma+z)}X_t(u)} p(u)\,du$$

$$\leq \frac{\alpha^2}{4(\gamma+z)^2} K(1) \quad , \quad (67)$$

$$\leq \frac{\alpha^2}{4\gamma^2} K(1)$$

which proves (33). The uniform continuity of $\frac{\partial H}{\partial z}$ follows from the last line of (66). The decreasing property of $H$ for $z \geq 0$ is immediate from its functional form.

□

### *Proof of Proposition 2*

This proof is an adaptation of the classical method of doubling variables using **Lemma 2**. Indeed, we require the following inequality for any viscosity super-solution $\overline{\Phi}$ and any sub-solution $\underline{\Phi}$:

$$\underline{\Phi} \leq \overline{\Phi} \text{ on } \overline{\Omega}. \quad (68)$$

The inequality is clearly satisfied along the boundaries $(0,T)\times\{0\}$ and $\{T\}\times[0,+\infty)$. Subsequently, we prove the inequality (68) at each point in $\Omega$ using the auxiliary function (e.g., Proof of Theorem 7.5 of [54]).

□

*Proof of Proposition 3*

The proof is owing to the fact that $\tilde{H}_\mu \to \tilde{H}_0$ as $\mu \to +0$ locally uniformly in $[0,T] \times \mathbb{R}$, combined with the continuity of the function $-R(n)nz + \tilde{H}(t,z)$ of $(t,n,z) \in [0,T] \times [0,+\infty) \times \mathbb{R}$ that is the full Hamiltonian of the HJBI equation, and he continuity of $\tilde{H}_0(t,z)$ with respect to all $(t,z) \in [0,T] \times \mathbb{R}$. The proposition then follows from the definition of viscosity solutions (e.g., Lemma 3.2 of [41]; Section 6 of [36]).

□

*Proof of Proposition 4*

We use a method of induction. From (41), we obtain $\Phi_{I,j} \geq 0$ ($0 \leq j \leq J$). Assume that we have $\Phi_{i+1,j} \geq 0$ ($0 \leq j \leq J$) for some $0 \leq i \leq I-1$. Using (44) and **Lemma 2**, we obtain

$$\begin{aligned}
\Phi_{i,j} &= \Phi_{i+1,j} - R(j\Delta n) j\Delta n \Delta t \frac{\Phi_{i+1,j} - \Phi_{i+1,j-1}}{\Delta n} + \tilde{H}\left(i\Delta t, \frac{\Phi_{i+1,j} - \Phi_{i+1,j-1}}{\Delta n}\right)\Delta t \\
&= \left(1 - R(j\Delta n) j\Delta n \Delta t \frac{1}{\Delta n}\right)\Phi_{i+1,j} + R(j\Delta n) j\Delta n \Delta t \frac{1}{\Delta n}\Phi_{i+1,j-1} + \tilde{H}\left(i\Delta t, \frac{\Phi_{i+1,j} - \Phi_{i+1,j-1}}{\Delta n}\right)\Delta t \\
&\geq \left(1 - R(j\Delta n) j\Delta t\right)\Phi_{i+1,j} + \tilde{H}\left(i\Delta t, \frac{\Phi_{i+1,j} - \Phi_{i+1,j-1}}{\Delta n}\right)\Delta t \\
&\geq \left(1 - R(j\Delta n) j\Delta t\right)\Phi_{i+1,j} + \tilde{H}\left(i\Delta t, \frac{\Phi_{i+1,j}}{\Delta n}\right)\Delta t \\
&\geq \left(1 - R(j\Delta n) j\Delta t\right)\Phi_{i+1,j} + \tilde{H}(i\Delta t, 0)\Delta t - \frac{\alpha^2}{4\gamma^2}K(1)\frac{\Phi_{i+1,j}}{\Delta n}\Delta t \\
&\geq \left(1 - R(j\Delta n) j\Delta t\right)\Phi_{i+1,j} - \frac{\alpha^2}{4\gamma^2}K(1)\frac{\Phi_{i+1,j}}{\Delta n}\Delta t \\
&= \left(1 - R(j\Delta n) j\Delta t - \frac{\alpha^2}{4\gamma^2}K(1)\frac{\Delta t}{\Delta n}\right)\Phi_{i+1,j}
\end{aligned} \qquad (69)$$

where we used $\Phi_{i+1,j} \geq 0$ ($0 \leq j \leq J$) and the Taylor expansion that is satisfied by some $y \geq 0$ (here, $\frac{\partial \tilde{H}}{\partial z}$ is a right derivative)

$$\begin{aligned}
\tilde{H}\left(i\Delta t, \frac{\Phi_{i+1,j}}{\Delta n}\right) &= \tilde{H}(i\Delta t, 0) + \frac{\partial \tilde{H}}{\partial z}(i\Delta t, y)\frac{\Phi_{i+1,j}}{\Delta n} \\
&\geq \tilde{H}(i\Delta t, 0) - \frac{\alpha^2}{4\gamma^2}K(1)\frac{\Phi_{i+1,j}}{\Delta n}
\end{aligned} \qquad (70)$$

and

$$\tilde{H}(i\Delta t, 0) = -\mu \ln\left(\int_0^1 e^{-\frac{\alpha^2}{4\mu\gamma}X_t(u)} p(u) \mathrm{d}u\right)$$

$$\geq -\mu \ln\left(\int_0^1 p(u) \mathrm{d}u\right) \tag{71}$$

$$= 0$$

Therefore, from the last line of (69), we obtain $\Phi_{i,j} \geq 0$ ($0 \leq j \leq J$) if

$$\Delta t \leq \Delta n \left(R(j\Delta n) j\Delta n + \frac{\alpha^2}{4\gamma^2} K(1)\right)^{-1} \quad \text{for all } 0 \leq j \leq J, \tag{72}$$

which is satisfied if

$$\Delta t \leq \Delta n \left(R(M) M + \frac{\alpha^2}{4\gamma^2} K(1)\right)^{-1}. \tag{73}$$

Consequently, we arrive at the desired result (45) by induction.

□

*Proof of Proposition 5*

We again use a method of induction. Using (41), we obtain $\Phi_{I,j} \leq h(J\Delta n) = h(M)$ ($0 \leq j \leq J$). Assume that we have $\Phi_{i+1,j} \leq \Xi$ ($0 \leq j \leq J$) for some $0 \leq i \leq I-1$ with some positive constant $\Xi$ that is independent from $\Delta t, \Delta n$. From (44), **Lemma 2**, and (69), we obtain

$$\begin{aligned}
\Phi_{i,j} &= \Phi_{i+1,j} - R(j\Delta n) j\Delta n \Delta t \frac{\Phi_{i+1,j} - \Phi_{i+1,j-1}}{\Delta n} + \tilde{H}\left(i\Delta t, \frac{\Phi_{i+1,j} - \Phi_{i+1,j-1}}{\Delta n}\right)\Delta t \\
&= \left(1 - R(j\Delta n) j\Delta n \Delta t \frac{1}{\Delta n}\right)\Phi_{i+1,j} + R(j\Delta n) j\Delta n \Delta t \frac{1}{\Delta n}\Phi_{i+1,j-1} + \tilde{H}\left(i\Delta t, \frac{\Phi_{i+1,j} - \Phi_{i+1,j-1}}{\Delta n}\right)\Delta t \\
&\leq \left(1 - R(j\Delta n) j\Delta t\right)\Phi_{i+1,j} + R(j\Delta n) j\Delta n \Delta t \frac{1}{\Delta n}\Xi + \tilde{H}\left(i\Delta t, \frac{\Phi_{i+1,j} - \Phi_{i+1,j-1}}{\Delta n}\right)\Delta t \\
&\leq \left(1 - R(j\Delta n) j\Delta t\right)\Phi_{i+1,j} + R(j\Delta n) j\Delta n \Delta t \frac{1}{\Delta n}\Xi + \tilde{H}\left(i\Delta t, \frac{\Phi_{i+1,j} - \Xi}{\Delta n}\right)\Delta t \\
&\leq \left(1 - R(j\Delta n) j\Delta t\right)\Phi_{i+1,j} + R(j\Delta n) j\Delta n \Delta t \frac{1}{\Delta n}\Xi + \tilde{H}(i\Delta t, 0)\Delta t + \frac{\alpha^2}{4\gamma^2}K(1)\frac{\Xi - \Phi_{i+1,j}}{\Delta n}\Delta t \\
&\leq \left(1 - R(j\Delta n) j\Delta t - \frac{\alpha^2}{4\gamma^2}K(1)\frac{\Delta t}{\Delta n}\right)\Phi_{i+1,j} + \left(R(j\Delta n) j\Delta n \Delta t \frac{1}{\Delta n} + \frac{\alpha^2}{4\gamma^2}K(1)\frac{\Delta t}{\Delta n}\right)\Xi \\
&\quad + H(i\Delta t, 0)\Delta t \\
&= \left(1 - R(j\Delta n) j\Delta t - \frac{\alpha^2}{4\gamma^2}K(1)\frac{\Delta t}{\Delta n}\right)\Xi + \left(R(j\Delta n) j\Delta n \Delta t \frac{1}{\Delta n} + \frac{\alpha^2}{4\gamma^2}K(1)\frac{\Delta t}{\Delta n}\right)\Xi \\
&\quad + H(i\Delta t, 0)\Delta t \\
&\leq \Xi + H(i\Delta t, 0)\Delta t
\end{aligned} \tag{74}$$

where we used $\Phi_{i+1,j} \leq \Xi$ ($0 \leq j \leq J$) and the Taylor expansion that is satisfied by some $y \geq 0$ (here,

$\dfrac{\partial \tilde{H}}{\partial z}$ is a right derivative):

$$\tilde{H}\left(i\Delta t, \dfrac{\Phi_{i+1,j} - \Xi}{\Delta n}\right) = \tilde{H}(i\Delta t, 0) + \dfrac{\partial \tilde{H}}{\partial z}(i\Delta t, y)\dfrac{\Phi_{i+1,j} - \Xi}{\Delta n}$$

$$= \tilde{H}(i\Delta t, 0) - \dfrac{\partial \tilde{H}}{\partial z}(i\Delta t, y)\dfrac{\Xi - \Phi_{i+1,j}}{\Delta n} \ . \tag{75}$$

$$\leq \tilde{H}(i\Delta t, 0) + \dfrac{\alpha^2}{4\gamma^2} K(1)\dfrac{\Xi - \Phi_{i+1,j}}{\Delta n}$$

By repeating (74) combined with $\Phi_{I,j} \leq h(M)$ ($0 \leq j \leq J$),

$$\Phi_{0,j} \leq h(M) + \Delta t \sum_{i=0}^{I-1} H(i\Delta t, 0), \ \ 0 \leq j \leq J \ . \tag{76}$$

We further obtain

$$\Delta t \sum_{i=0}^{I-1} H(i\Delta t, 0) = \Delta t \sum_{i=0}^{I-1}\left(-\mu \ln\left(\int_0^1 e^{-\frac{\alpha^2}{4\mu\gamma}X_t(u)} p(u)\mathrm{d}u\right)\right)$$

$$\leq \Delta t \sum_{i=0}^{I-1}\left(-\mu \ln\left(\int_0^1 e^{-\frac{\alpha^2}{4\mu\gamma}K(1)} p(u)\mathrm{d}u\right)\right)$$

$$= \Delta t \sum_{i=0}^{I-1}\left(-\mu\left(-\dfrac{\alpha^2}{4\mu\gamma}K(1)\right)\right) \ . \tag{77}$$

$$= I\Delta t \dfrac{\alpha^2}{4\gamma}K(1)$$

$$= \dfrac{\alpha^2 T}{4\gamma}K(1)$$

We obtain the desired result directly from (76) and (77); namely, we should set $\Xi$ larger than the right-hand side of (76).

□

## *Proof of Proposition 6*

The first statement is immediately verified by the partial differentiations

$$\frac{\partial}{\partial \Phi_{i+1,j}}\left(\Phi_{i+1,j} - R(j\Delta n)j\Delta n\Delta t \frac{\Phi_{i+1,j} - \Phi_{i+1,j-1}}{\Delta n} + \tilde{H}\left(i\Delta t, \frac{\Phi_{i+1,j} - \Phi_{i+1,j-1}}{\Delta n}\right)\Delta t\right)$$

$$= 1 - R(j\Delta n)j\Delta n\Delta t \frac{1}{\Delta n} - \frac{\partial}{\partial \Phi_{i+1,j}} \tilde{H}\left(i\Delta t, \frac{\Phi_{i+1,j} - \Phi_{i+1,j-1}}{\Delta n}\right)\Delta t$$

$$\geq 1 - R(j\Delta n)j\Delta n\Delta t \frac{1}{\Delta n} - \frac{\alpha^2}{4\gamma^2} K(1) \frac{\Delta t}{\Delta n} \tag{78}$$

$$= \frac{1}{\Delta n}\left(\Delta n - \left(R(j\Delta n)j\Delta n + \frac{\alpha^2}{4\gamma^2} K(1)\right)\Delta t\right)$$

$$\geq \frac{1}{\Delta n}\left(\Delta n - \left(R(M)M + \frac{\alpha^2}{4\gamma^2} K(1)\right)\Delta t\right)$$

$$\geq 0$$

$$\frac{\partial}{\partial \Phi_{i+1,j-1}}\left(\Phi_{i+1,j} - R(j\Delta n)j\Delta n\Delta t \frac{\Phi_{i+1,j} - \Phi_{i+1,j-1}}{\Delta n} + \tilde{H}\left(i\Delta t, \frac{\Phi_{i+1,j} - \Phi_{i+1,j-1}}{\Delta n}\right)\Delta t\right)$$

$$= R(j\Delta n)j\Delta n\Delta t \frac{1}{\Delta n} + \frac{\partial}{\partial \Phi_{i+1,j-1}} \tilde{H}\left(i\Delta t, \frac{\Phi_{i+1,j} - \Phi_{i+1,j-1}}{\Delta n}\right)\Delta t \tag{79}$$

$$\geq R(j\Delta n)j\Delta n\Delta t \frac{1}{\Delta n} + \frac{\Phi_{i+1,j-1}}{\Delta n}\Delta t\left\{-\frac{\partial \tilde{H}}{\partial z}\left(i\Delta t, \frac{\Phi_{i+1,j} - \Phi_{i+1,j-1}}{\Delta n}\right)\right\}$$

$$\geq 0$$

where we used $\Phi_{i+1,j-1} \geq 0$ and **Lemma 2**. Note that $\frac{\partial \tilde{H}}{\partial z}$ is understood as the right derivative in the inequalities above. The second statement concerning the increasing property follows from the proof of Lemma 5.2, 3(a) of [55].

□

**Appendix B: Few additional computational results**

We discussed influences of the uncertainty aversion on the growth curve in **Figure 7**. In this appendix, we firstly demonstrate how the value function and the optimal harvesting policy look like when there is no uncertainty ($\mu \to +\infty$). We assume the uniform distribution $p$ here. **Figures B1(a)** and **(b)** show the computed value function and optimal harvesting policy for 2021 with sustainability concern when there is no uncertainty, respectively. The computational conditions other than $\mu$ (1/day) are the same with that for the model with sustainability concern in 2021 (e.g., **Figures 4(a)** and **6(a)**), and we set $\mu = 10^{10}$ here that is a sufficiently large value to mimic the case $\mu \to +\infty$. Comparing **Figure 4(a)** with **Figure B1(a)** reveals that the value function is smaller for the larger uncertainty aversion, which is considered due to the underestimated body weight in the former. Comparing **Figure 6(a)** with **Figure B1(b)** reveals that the maximum harvesting intensity becomes larger for the smaller uncertainty aversion.

Secondly, we analyze the worst-case distorted distribution, which is $p$ multiplied by $\hat{\phi}$ in (26), when $p$ it is the beta type. For the computation here, the integration with respect to $u$ that is involved in $\tilde{H}$ is carried out by a midpoint rule with uniformly distributed at 300 points to create figures with a higher resolution. **Figures B2(a)-(c)** show the distorted and non-distorted beta distributions for the right-skewed case $(a,b)=(2,5)$, non-skewed case $(a,b)=(2,2)$, and left-skewed case $(a,b)=(5,2)$. We have computed the distorted distribution at the time $t=61$ (day) with $n=5$ as demonstrative examples. These figures demonstrate that increasing the uncertainty aversion shifts the mode of the original beta distributions toward left with the increased height. The case $\mu=0.1$ is close to the no-uncertainty case and the position of the maximizer of the distributions remain almost the same, whereas the case $\mu=0.01$ leads to the distortion more clearly visible. The distorted distributions are still smooth and bounded. This kind of visualization analysis may help model identifications in applied studies.

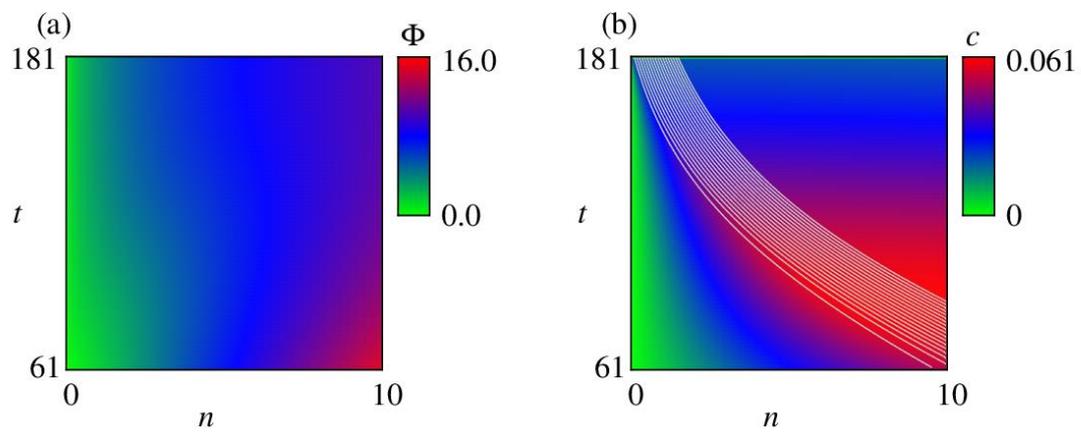

**Figure B1.** Computed results for the year 2021 with sustainability concern: **(a)** value function and **(b)** optimal harvesting policy and the population trajectories.

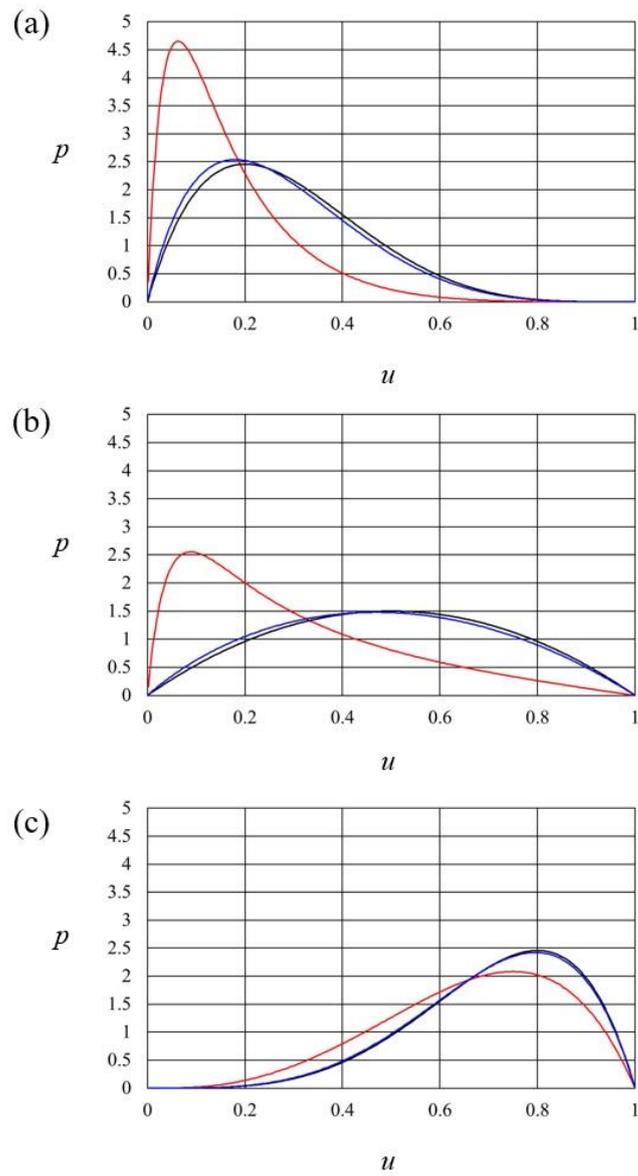

**Figures B2.** The distorted and non-distorted beta distributions, denoted as $p$, for the right-skewed case (a) $(a,b)=(2,5)$, (b) non-skewed case $(a,b)=(2,2)$, and (c) left-skewed case $(a,b)=(5,2)$. Each color represents $\mu=0.01$ (red), $\mu=0.1$ (blue), and $\mu=+\infty$ (black).